\documentclass[final,3p,times,sort&compress,doublespacing,fleqn]{elsarticle}

\journal{arXiv}

\usepackage{multirow,multicol,setspace}

\usepackage{amsthm,amssymb,array,float,latexsym,amsmath,amssymb,amsbsy,eurosym,theorem,amsfonts,mathrsfs,bbm,verbatim,dsfont}
\usepackage{booktabs}
\usepackage{caption}
\usepackage{subcaption}

\usepackage{rotating}
\usepackage{tabularx}
\usepackage[table]{xcolor}
\usepackage[]{graphicx}
\usepackage{color}

\usepackage[most]{tcolorbox}
\usepackage{titlesec}
\usepackage{algorithm}
\usepackage{algpseudocode}

\usepackage{scalerel,stackengine}

\titleformat{\paragraph}
{\normalfont\normalsize}{\theparagraph}{1em}{}
\titlespacing*{\paragraph}
{0pt}{3.25ex plus 1ex minus .2ex}{1.5ex plus .2ex}

\renewcommand{\b}[1]{\boldsymbol{#1}} 
\newcommand{\cc}[1]{\mathcal{#1}}

\renewcommand{\v}[1]{\mathbbm{#1}} 
\renewcommand{\v}[1]{\text{\usefont{U}{bbm}{m}{n}#1}} 
\newcommand{\h}[1]{\widehat{#1}}
\renewcommand{\o}[1]{\overline{#1}}

\renewcommand{\t}[1]{\widetilde{#1}}

\newcommand{\lp}{\left(}
\newcommand{\rp}{\right)}

\newcommand{\lb}{\left[}
\newcommand{\rb}{\right]}

\newcommand{\lc}{\left\{}
\newcommand{\rc}{\right\}}

\newcommand{\p}{\partial}

\newcommand{\f}{\displaystyle\frac}

\newcommand{\bnull}{\b0}

\renewcommand{\d}{\mbox{d}}

\newcommand{\curl}{\mbox{curl}}

\newcommand{\dyad}{\otimes}

\definecolor{gray}{rgb}{0.75, 0.75, 0.75}
\definecolor{yellow}{rgb}{1, 0.7, 0.2}
\definecolor{green}{rgb}{0.3, 0.9, 0.3}
\definecolor{brown}{rgb}{0.6, 0.3, 0.2}
\definecolor{magenta}{rgb}{0.9, 0.1, 0.9}
\definecolor{light}{rgb}{1, 0.7, 0.7}

\graphicspath{{figs/}}

\newcolumntype{Y}{>{\centering\arraybackslash}X}

\begin{document}

\begin{frontmatter}

\title{A micromorphic-based artificial diffusion method for stabilized finite element approximation of convection-diffusion problems}

\author[erlangen]{Soheil~Firooz\corref{cor}}
\ead{soheil.firooz@fau.de}
\author[Africa]{B.~Daya~Reddy}
\author[erlangen,glasgow]{Paul~Steinmann}
\address[erlangen]{Institute of Applied Mechanics, Friedrich-Alexander-Universit{\"a}t Erlangen--N{\"u}rnberg , Egerland Str. 5, 91058 Erlangen, Germany}
\address[Africa]{Department of Mathematics and Applied Mathematics and Centre for Research in Computational and Applied Mechanics, University of Cape Town, Rondebosch, 7701, South Africa}
\address[glasgow]{Glasgow Computational Engineering Center, James Watt School of Engineering, University of Glasgow, Glasgow G12 8QQ, United Kingdom}
\cortext[cor]{Corresponding author.}

\begin{abstract}
We present a novel artificial diffusion method to circumvent the instabilities associated with the standard finite element approximation of convection-diffusion equations.
Motivated by the micromorphic approach, we introduce an auxiliary variable, which is related to the gradient of the field of interest, and which leads to a coupled problem.
Conditions for well-posedness of the resulting formulation are established.
We carry out a comprehensive numerical study to compare the proposed methodology against some well-established approaches in one- and two-dimensional settings.
The proposed method outperforms established approaches in general in approximating accurately the solutions to pertinent and challenging problems.
\end{abstract}

\begin{keyword}
Convection-diffusion problems, Artificial diffusion, Convection-dominated problems, Stability
\end{keyword}

\end{frontmatter}

\section{Introduction}\label{sec:intro}

\subsection{Convection-diffusion}\label{sec:conv-diff}
\noindent
Convection-diffusion problems are fundamental in many areas including mechanics, chemistry and biology, and include the modelling of phenomena such as diffusion of solutes in liquids, heat transfer in thermal systems, and turbulent flow over aircraft wings.
If the diffusive effects are significantly smaller than the convective effects, convection-diffusion problems are commonly referred to as convection-dominated~\cite{Christie1976,Heinrich1977,Johnson1984}.
In such scenarios, the transport of quantities such as heat, mass, or momentum is primarily governed by their convection rather than their diffusion.
This characteristic can result in the formation of sharp gradients in the solution, posing challenges for numerical simulation of the problem.

\subsection{Finite element technology}\label{sec:FE}
\noindent
It is well-established that the standard Galerkin finite element method (FEM) often encounters instabilities when dealing with convection-dominated problems.
Such instabilities lead to oscillations or unwanted diffusion, compromising the accuracy and reliability of simulations~\cite{Gresho1981,Brooks1982,Johnson1984}.
This situation has motivated the development of alternative methodologies that avoid spurious oscillations regardless of the degree of mesh refinement.

For convection-dominated problems it was originally observed that oscillation-free solutions are obtainable via upwinding of the convective term~\cite{Christie1976,Heinrich1977}.
The upwinded convective term can be constructed by simply adding artificial diffusion~\cite{Brooks1982}.
While this method helps to obtain stable solutions, it can introduce excessive smoothing (smearing) of sharp gradients or boundary layers, depending on how much artificial diffusion is added.
Hughes~\cite{Hughes1978} proposed a reduced quadrature method as an alternative formulation to the classical upwinding.
Therein the numerical quadrature rule for the convection term is modified to achieve the upwind effect.
Despite precluding oscillations, this method still suffers from a lack of accuracy due to over-smoothing of the solution.

The streamline upwind (SU) and streamline upwind Petrov--Galerkin (SUPG) methods~\cite{Brooks1982,Hughes1982b} incorporate the streamline direction into the test function.
Via a stabilization parameter, the fluxes may be adjusted to effectively capture the convective effects.
Consequently, both methods effectively damp out spurious oscillations while maintaining stability.
Although the original SU and SUPG methods proved to be stable and accurate, in problems with non-smooth solutions such as those involving boundary layers or sharp gradients, spurious oscillations, undershoots, and overshoots could still emerge near regions with steep gradients~\cite{Burman2010}, specially if the flow direction is skew to the mesh.
Hughes et al.~\cite{Hughes1986,Hughes1986d} addressed this issue by introducing further a discontinuity-capturing term to the weighting function, which was analogous to the streamline term, but targeting the direction of the solution gradient rather than the flow direction.
See~\cite{Knopp2002,DeSampaio2001,Codina1993a} for further studies on discontinuity-capturing and shock-capturing methods.

Later, Hughes and his coworkers~\cite{Hughes1995,Hughes1998,Buffa2006,John2008a} introduced the variational multiscale (VMS) method which decomposes the solution into coarse (resolved) and fine (unresolved) scales and then systematically models the effect of the fine scales on the resolved ones.
This approach is capable recovering classical stabilization schemes such as SUPG method under certain conditions.
In the context of convection-diffusion problems, VMS enables better control of numerical oscillations and enhanced accuracy of fluxes, particularly for convection-dominated cases.
Followed by the work of Layton~\cite{Layton2002}, John et al.~\cite{John2006a} introduced a two-level variational multiscale method where they employed a fine mesh finite element space to approximate the concentration and a coarse-mesh discontinuous finite element space to represent the large-scale components of the flux in the two-scale discretization, see also~\cite{Chen2018a,Du2015}.
Another extension of the VMS method is the local projection stabilization (LPS) method where certain components of the solution (usually the gradient or convective terms) are locally projected onto a discontinuous or coarser finite element space, and then the difference between the original quantity and its projection is penalized~\cite{Knobloch2009,Cibik2011,Matthies2009}.

Aside from the SUPG and VMS methods, the discontinuous Galerkin finite element method (DG-FEM) offers a robust approach to overcoming instabilities in convection-diffusion equations~\cite{Johnson1986,Cockburn1999c}.
Unlike continuous Galerkin methods, DG-FEM allows for discontinuities between elements, providing greater flexibility and control over the numerical solution.

Further examples of stabilizing methodologies for convection-diffusion equations include the Galerkin least squares (GLS) method~\cite{Hughes1989a,Pironneau1992}, the symmetric stabilization method~\cite{Burman2009,Guermond1999,Codina2000a}, the residual-free bubbles (RFB) method~\cite{Franca1998b,Brezzi1999b}, the generalized Taylor--Galerkin method~\cite{Donea1984}, the mixed-order finite element method~\cite{Taylor1973,Case2011} and the continuous interior penalty method.

\subsection{Micromorphic-based artificial diffusion}\label{sec:MMAD-intro}
\noindent
In this contribution, motivated by the Mean Zero Artificial Diffusion (MZAD) method introduced in~\cite{Firooz2024}, we develop a MicroMorphic-based Artificial Diffusion (MMAD) method for stable and accurate finite element approximation of convection-diffusion problems.
The mean zero artificial diffusion method is equivalent to projected artificial diffusion approach, with the artificial diffusivity parameter having to be large enough to guarantee the stability and well-posedness of the solution.
As will be elucidated later, this requirement can lead to solutions with excessive wiggles in the vicinity of very steep gradients..
This issue gets worse for higher dimensional problems since the artificial diffusivity parameter does not take the flow direction into account.

The micromorphic-based artificial diffusion (MMAD) method proposed here is, in principle, a gradient-enhanced version of the mean zero artificial diffusion method in which the artificial diffusion takes the flow direction into account and is capable of capturing sharp gradients.
Motivated by the micromorphic approach~\cite{Forest2009,Forest2016,Forest2020}, we modify the convection-diffusion problem by introducing a new micromorphic-type variable and carrying out enhancement in terms of the variable and its gradient.
This equips the conventional convection-diffusion problem with a fictitious potential-like function, which is analogous to that of the micromorphic approach, but differs in the underlying principle.

Equal-order conforming approximations are adopted in the finite element formulation and the resulting discrete problem is shown to be well-posed, with convergence independent of the diffusivity.
To evaluate the performance of our proposed methodology, we carry out a comprehensive numerical study using a number of benchmark examples and compare our approach with other established methods.
Our proposed methodology shows improved performance compared to the VMS, SU, SUPG and MZAD methods.

Table~\ref{tab:Abb} provides a list of abbreviations used throughout this manuscript.

\begin{table}
\centering
\begin{tabular}{|p{0.15\textwidth}  p{0.6\textwidth}|}
\toprule
\multicolumn{2}{|p{0.5\textwidth}|}{\textbf{Abbreviations}}\\
\midrule
SUPG     & streamline upwind Petrov--Galerkin method \\
MZAD     & mean zero artificial diffusion method \\
MMAD     & micromorphic-based artificial diffusion method \\
\bottomrule
\end{tabular}
\caption
{ 
List of frequently used abbreviations used throughout this manuscript.
}
\label{tab:Abb}
\end{table}

\subsection{Organization of the manuscript}\label{sec:org}
\noindent
The remainder of this manuscript is organized as follows: 
Section~\ref{sec:governing equations} introduces the problem definition and presents the governing equations together with recapitulating well-established methodologies to overcome instabilities associated with convection-diffusion problems.
The micromorphic-based artificial diffusion method is elaborated in Section~\ref{sec:MMAD}. 
Section~\ref{sec:Discretization} presents details of the well-posedness and convergence of our approach.
Section~\ref{sec:result} illustrates the proposed MMAD method through a set of numerical examples together with comparison with alternative approaches. 
Section~\ref{sec:conc} summarizes the work and presents an outlook of further work.


\section{Convection-diffusion problems}\label{sec:governing equations}

\subsection{Problem statement}\label{sec:problem}
\noindent
Let $\cc{B}$ be a bounded domain with smooth boundary $\p \cc{B}$.
The general \emph{dimensionless} form of the convection-diffusion equation for a scalar quantity $\varphi$ reads
\begin{equation}
\begin{aligned}
&\f{\p\varphi}{\p t} +  \b{u} \cdot \nabla \varphi  -  k \nabla \cdot \lp\nabla \varphi\rp   = F  \quad \text{in}\,\, \cc{B} \,,  \\
&\varphi = \varphi_{\text{p}}  \quad \text{on}\,\, \p\cc{B}_{\text{D}} \, 
\qquad \text{and} \qquad 
\nabla \varphi \cdot \b{n} = \f{\p \varphi}{\p \b{n}} = \b{t}_{\text{p}} \quad \text{on}\,\, \p\cc{B}_{\text{N}} \,, \\
&\varphi(x,0)=\varphi_{0}(x)\quad \text{on}\,\, \cc{B}\,, \\
&k=\f{1}{\text{Pe}}=\f{D}{L|\b{u}|} \,,
\end{aligned}
\label{eq:strong-1}
\end{equation}
with $\b{u}$ the divergence-free velocity field, $D$ the diffusion coefficient, $L$ the characteristic length of the domain, Pe the Peclet number, $\b{n}$ the unit vector normal to the boundary $\p \cc{B}$, $\varphi_{\text{p}}$ the prescribed value of $\varphi$, $\varphi_{0}$ the initial condition and $\b{t}_{\text{p}}$ the prescribed flux.
$\p\cc{B}_{\text{D}}$ and $\p\cc{B}_{\text{N}}$ denote those parts of the boundary $\p\cc{B}$ on which Dirichlet and Neumann boundary conditions are prescribed, respectively with $\p\cc{B}_{\text{D}}\cap\p\cc{B}_{\text{N}}=\emptyset$ and $\p\cc{B}_{\text{D}}\cup\p\cc{B}_{\text{N}}=\p\cc{B}$.
The first term in Eq.~\eqref{eq:strong-1} accounts for the time dependence of the problem, the second term accounts for convection, the third term accounts for diffusion and the term on the right is a source term.
The Peclet number characterizes the relative significance of convective and diffusive terms in Eq.~\eqref{eq:strong-1}.
When $\text{Pe} \ll 1$, diffusion dominates over convection and the solution exhibits smooth profiles.
When $\text{Pe} \gg 1$, convection dominates over diffusion and the solution is more likely to exhibit sharp fronts and gradients.

Taking the inner product of Eq.~\eqref{eq:strong-1} with a test function $\delta \varphi$ that satisfies the homogeneous Dirichlet boundary condition, and using the Neumann boundary condition, the weak form of Eq.~\eqref{eq:strong-1} reads
\begin{equation}
\int_{\cc{B}}   \f{\p\varphi}{\p t}\, \delta \varphi \, \d v + \int_{\cc{B}}  \b{u} \cdot \nabla \varphi \, \delta \varphi \, \d v + \int_{\cc{B}} k \nabla \varphi \cdot  \nabla \delta \varphi \, \d v  = \int_{\cc{B}} F \,  \delta \varphi \, \d v + \int_{\p\cc{B}_{\text{N}}} \lb \nabla \varphi \cdot \b{n} \rb \, \delta \varphi   \, \d \b{a} \,.
\label{eq:weak-tmp}
\end{equation}

\begin{figure}[b!]
\centering
\includegraphics[width=1.0\textwidth]{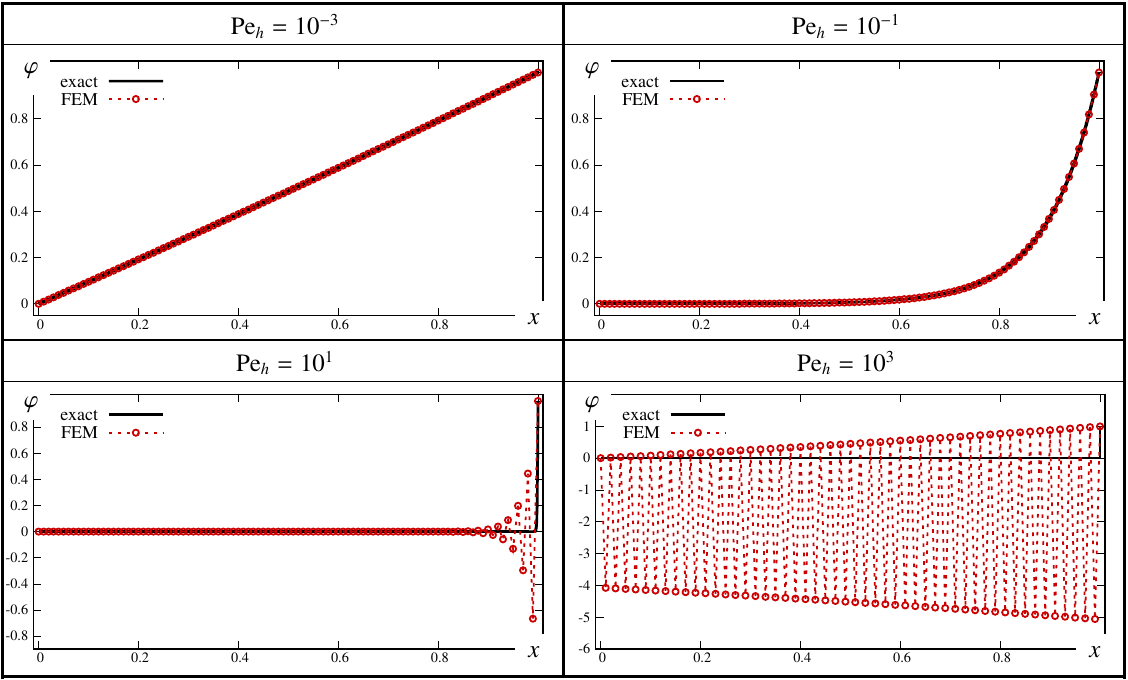}
\caption{ 
Standard Galerkin FEM solution versus the exact solution in a 1D convection–diffusion problem at different Peclet numbers.
}
\label{fig:FEM}
\end{figure}

\subsection{Stabilizing approaches for convection-diffusion problems}\label{sec:instability}
\noindent
In this work, we will consider various finite element approximations of Eq.~\eqref{eq:weak-tmp}.
With the mesh size denoted by $h$, it is well-known that if $h>2D/|\b{u}|$, the standard finite element solution displays spurious oscillations~\cite{Hughes1978}, i.e. the method fails to accurately capture the sharp gradients that are associated with high Peclet numbers.
It proves to be more relevant to define the element Peclet number as ${\text{Pe}_{h}}=h|\b{u}|/D$ as it accounts for local mesh resolution and provides a more accurate characterization of the balance between convection and diffusion at the element level.
For illustration and as a reference, Fig.~\ref{fig:FEM} shows the one-dimensional finite element solution for a steady-state convection-diffusion problem with no source term.
A domain of unit length is discretized into $100$ linear finite elements and is subject to Dirichlet boundary conditions $\varphi=0$ at $x=0$ and $\varphi=1$ at $x=1$.
The solid black line shows the exact solution and the dashed red line (with points on top) shows the discrete solution.
It is observed that for small Peclet numbers the discrete solution accurately captures the exact solution.
However, as the Peclet number increases and the problem becomes more convection-dominated, instabilities in the form of wiggles appear in the finite element solution.

The classical artificial diffusion method adds artificial diffusivity to the diffusion term in order to eliminate the oscillations in the solution.
That is, the strong form is modified to
\begin{equation}
\int_{\cc{B}}   \f{\p\varphi}{\p t}\, \delta \varphi \, \d v + \int_{\cc{B}}  \b{u} \cdot \nabla \varphi \, \delta \varphi \, \d v + \int_{\cc{B}} \lb k + \o{k} \rb \nabla \varphi \cdot \nabla \delta \varphi \, \d v = \int_{\cc{B}} F \,  \delta \varphi \, \d v + \int_{\p\cc{B}_{\text{N}}}  \b{t}_{\text{p}} \, \delta \varphi   \, \d \b{a} \,,
\label{eq:weak-AD}
\end{equation}
with $\o{k}$ the artificial diffusion coefficient.

The streamline upwinding (SU) method~\cite{Brooks1982,Brooks1981} is similar to the classical artificial diffusion method where an artificial diffusion coefficient $\o{k}$ is chosen by taking the flow direction into account.
In this approach, $\o{k}$ is a scalar in a one-dimensional problem whereas in a multi-dimensional case, it is a second order tensor referred to as the artificial diffusivity tensor.
The SU method is implemented by modifying the test function associated with the convective term as 
\begin{equation}
\int_{\cc{B}}   \f{\p\varphi}{\p t}\, \delta \varphi \, \d v + \int_{\cc{B}}  \b{u} \cdot \nabla \varphi \, \o{\delta \varphi} \, \d v + \int_{\cc{B}} k \nabla \varphi \cdot \nabla \delta \varphi \, \d v = \int_{\cc{B}} F \,  \delta \varphi \, \d v + \int_{\p\cc{B}_{\text{N}}}  \b{t}_{\text{p}} \, \delta \varphi   \, \d \b{a} \,,
\label{eq:weak-SUPG}
\end{equation}
with $\o{\delta \varphi}$ defined, for two-dimensional problems for example, as $\o{\delta \varphi} = \delta \varphi + \o{k} \, \b{u} \dyad \nabla \delta \varphi / |\b{u}|$ where
\begin{equation}
\begin{aligned}
&\o{k} = \sum_{i=1}^{\text{PD}} u_{i} h_{i} {\gamma}_{i} / 2\,,\qquad \text{with} \qquad
{\gamma}_{i} = \coth(\alpha_i)-1/\alpha_i\,,\qquad
\alpha_i=u_i h_i/2D\,,\qquad
u_i = \b{e}_i^T \cdot \b{u}\,,
\end{aligned}
\label{eq:KK0}
\end{equation}
with PD being the problem dimension, and $\b{e}_{i}$ being the unit vectors in the natural coordinate system used in FEM.
Application of $\o{\delta \varphi}$ to all the terms in Eq.~\eqref{eq:weak-tmp} yields the streamline upwind Petrov--Galerkin (SUPG)~\cite{Brooks1982} method with the weak form
\begin{equation}
\int_{\cc{B}}   \f{\p\varphi}{\p t}\, \o{\delta \varphi} \, \d v + \int_{\cc{B}}  \b{u} \cdot \nabla \varphi \, \o{\delta \varphi} \, \d v + \int_{\cc{B}} k \nabla \varphi \cdot \nabla \o{\delta \varphi} \, \d v = \int_{\cc{B}} F \,  \o{\delta \varphi} \, \d v + \int_{\p\cc{B}_{\text{N}}}  \b{t}_{\text{p}} \, \o{\delta \varphi}   \, \d \b{a}\,.
\label{eq:weak-SUPG-2}
\end{equation}


\section{The micromorphic-based artificial diffusion (MMAD) method}\label{sec:MMAD}

\noindent
Here we propose a novel micromorphic-based artificial diffusion (MMAD) approach for convection-diffusion or pure advection problems.

\subsection{Starting point: The mean zero artificial diffusion (MZAD) method~\cite{Firooz2024}}\label{sec:MZAD-test}
\noindent
First, to set the stage and for the sake of completeness, the mean zero artificial diffusion (MZAD) method is briefly introduced.
We introduce an additional variable $\b{g}$ to represent the gradient of $\varphi$, that is,
\begin{equation}
\b{g} := \nabla \varphi \Longrightarrow \nabla \varphi -\b{g} := \bnull \,.
\label{eq:g}
\end{equation}
For convenience we assume henceforth a homogeneous Dirichlet boundary condition on all of $\p B$.
We add the divergence of Eq.~\eqref{eq:g} to Eq.~\eqref{eq:strong-1} with a penalty parameter $p$, so that the modified problem becomes
\begin{subequations}
\begin{align}
& \f{\p\varphi}{\p t} +  \b{u} \cdot \nabla \varphi  -  k \Delta \varphi - p\nabla \!\cdot(\nabla \varphi - \b{g})   = F \,, 
\label{eq:MZAD-strong-1}\\
& \nabla \varphi -\b{g} = \bnull\,.\,
\label{eq:MZAD-strong-2}
\end{align}
\label{eq:MZAD-strong}
\end{subequations}
Following Eq.~\eqref{eq:weak-tmp} and multiplying Eq.~\eqref{eq:MZAD-strong-2} with a test function $\delta\b{g}$, the corresponding weak form reads
\begin{subequations}
\begin{align}
& \int_{\cc{B}}   \f{\p\varphi}{\p t}\, \delta \varphi \, \d v + \int_{\cc{B}}  \b{u} \cdot \nabla \varphi \, \delta \varphi \, \d v + \int_{\cc{B}} k \nabla \varphi \cdot \nabla \delta \varphi \, \d v  + \int_{\cc{B}} p \lb \nabla \varphi - \b{g} \rb\cdot \nabla \delta \varphi \,\d v = \int_{\cc{B}} F \,  \delta \varphi \, \d v\,,
\label{eq:MZAD-weak-1}\\
& \int_{\cc{B}}  \lb \nabla \varphi -\b{g} \rb \cdot \, \delta \b{g} \, \d v = 0\,.
\label{eq:MZAD-weak-2}
\end{align}
\label{eq:MZAD-weak}
\end{subequations}
We define $\Phi=\cc{H}_0^1(\cc{B})=\{\varphi \in L^{2}(\cc{B}), \, \p\varphi/\p x_i \in L^2(\cc{B}), \varphi=0 \; \text{on} \; \p\cc{B}   \}$ with the norm $||\varphi||_{\Phi}=||\nabla\varphi||_{L^2}=\lb\int_{\cc{B}} |\nabla \varphi |^2 \, \d v \rb^{1/2}$.
We also define $G=\{\b{g}\,|\, g_i \in L^2(\cc{B})  \}$ with the norm $||\b{g}||^2_{G} = \int_{\cc{B}} |\b{g}|^2 \, \d v$.
Then we seek $\varphi \in \Phi$ and $\b{g} \in G$ that satisfy Eq.~\eqref{eq:MZAD-weak} for all $\delta\varphi \in \Phi$ and $\delta\b{g} \in G$.
The resulting system of equations is consistent.
Indeed, Eq.~\eqref{eq:MZAD-weak-2} implies that the last term in Eq.~\eqref{eq:MZAD-weak-1} is zero.

For discrete approximations, we consider continuous piecewise linear finite element approximations on a mesh of triangles (2D) or tetrahedra (3D).
The discrete spaces are denoted by $\Phi^{h} \subset \Phi$ and $G^{h} \subset G$.
Then we seek $\varphi_{h} \in \Phi^{h}$ and $\b{g}_{h} \in G^{h}$ that satisfy
\begin{subequations}
\begin{align}
& \int_{\cc{B}}   \f{\p\varphi_{h}}{\p t}\, \delta \varphi_{h} \, \d v + \int_{\cc{B}}  \b{u} \cdot \nabla \varphi_{h} \, \delta \varphi_{h} \, \d v + \int_{\cc{B}} k \nabla \varphi_{h} \cdot \nabla \delta \varphi_{h} \, \d v  +\!\! \int_{\cc{B}}\!\! p \lb \nabla \varphi_{h} - \b{g}_{h} \rb\cdot \nabla \delta \varphi_{h} \,\d v = \int_{\cc{B}}\!\! F \,  \delta \varphi_{h} \, \d v\,,
\label{eq:MZAD-weak-discrete-1}\\
& \int_{\cc{B}}  \lb \nabla \varphi_{h} -\b{g}_{h} \rb \cdot \, \delta \b{g}_{h} \, \d v = 0\,,
\label{eq:MZAD-weak-discrete-2}
\end{align}
\label{eq:MZAD-weak-discrete}
\end{subequations}
for all $\delta\varphi_{h} \in \Phi^{h}$ and $\delta\b{g}_{h} \in G^{h}$.
In contrast to the continuous problem, Eq.~\eqref{eq:MZAD-weak-discrete-2} gives
\begin{equation}
\b{g}_{h} = \b{P} \cdot \nabla \varphi_{h}\,,
\label{eq:projection}
\end{equation}
where $\b{P}$ is the $L^2$-orthogonal projection of $\nabla \varphi_h$ onto $G^{h}$.
Substitution of Eq.~\eqref{eq:projection} in~\eqref{eq:MZAD-weak-discrete-1} yields
\begin{equation}
 \int_{\cc{B}}   \f{\p\varphi_{h}}{\p t}\, \delta \varphi_{h} \, \d v + \int_{\cc{B}}  \b{u} \cdot \nabla \varphi_{h} \, \delta \varphi_{h} \, \d v + \int_{\cc{B}} \lb k \nabla \varphi_{h} + p \lb \b{I} - \b{P} \rb \cdot \nabla \varphi_{h} \rb \cdot \nabla \delta \varphi_{h} \, \d v  = \!\!\int_{\cc{B}}\!\! F \,  \delta \varphi_{h} \, \d v\,,
\label{eq:MZAD-weak-discrete-1-updated}
\end{equation}
with $\b{I}$ being the second-order identity tensor.
Since $\lb \b{I} - \b{P} \rb$ is also a projection, we observe that the MZAD method is effectively a projected artificial diffusion formulation, see also~\cite{Layton2002,John2006a}.

The MZAD method proves to be stable for transient problems and is capable of capturing sharp gradients to some extent~\cite{Firooz2024}.
However, in the case of convection-dominated problems with extremely sharp gradients or extremely thin boundary layers, depending on the choice of the penalty parameter, it either suffers from lack of stability or over-smoothing of the solution.

\begin{figure}[b!]
\centering
\includegraphics[width=1.0\textwidth]{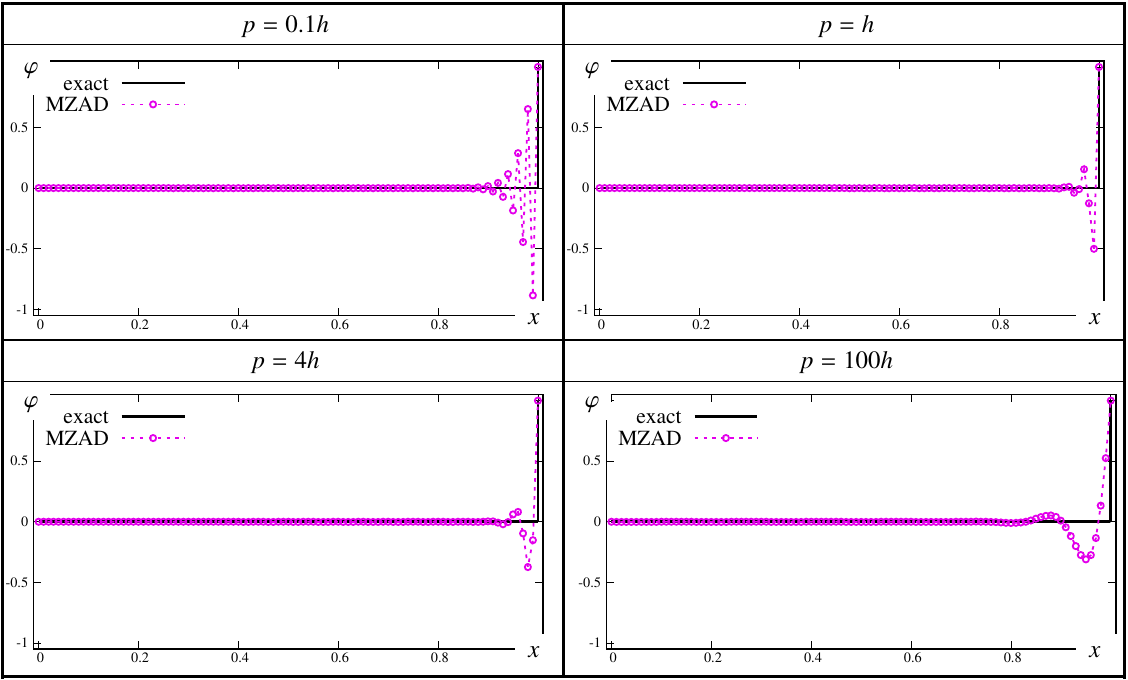}
\caption{ 
Mean zero artificial diffusion (MZAD) solution with different penalty parameters $p$ versus the exact solution in 1D convection-diffusion problems with $\text{Pe}_h=10^5$.
}
\label{fig:MZAD}
\end{figure}

Fig.~\ref{fig:MZAD} shows the solution of the MZAD method for the same convection-dominated problem as in Fig.~\ref{fig:FEM} with $\text{Pe}_{\text{h}}=10^5$.
We consider four different values of the penalty parameter $p$ as multiples of the mesh size $h$.
We observe that on the one hand small values of the penalty parameter yield slightly unstable solutions and on the other hand large values of the penalty parameter yield stable but over-smooth solutions.
Moreover, the parameter $p$ does not take the flow direction into account, which exacerbates its behavior in multi-dimensional problems.
We will elaborate in the following on how these two shortcomings are fully rectified by our proposed micromorphic-based artificial diffusion (MMAD) method which is, in particular, a generalized gradient-enhanced MZAD method.

\subsection{Improvement: The micromorphic-based artificial diffusion (MMAD) method}\label{sec:MZAD-test}
\noindent
To motivate the micromorphic-based artificial diffusion (MMAD) approach, consider the pure diffusion problem 
\begin{equation}
\int_{\cc{B}}  k \nabla \varphi \cdot \nabla \delta \varphi \, \d v = \int_{\cc{B}}\!\! F \,  \delta \varphi \, \d v\,.
\label{eq:pure-diffusion}
\end{equation}
This weak formulation corresponds to the minimization problem 
\begin{equation}
\underset{\varphi \in \Phi}{\text{min}} \; J_{0}(\varphi)\,, \quad  J_{0}(\varphi) = \f{k}{2} \int_{\cc{B}} |\nabla \varphi |^2 \, \d v + \int_{\cc{B}} F \varphi \, \d v\,.
\label{eq:minimization}
\end{equation}
Next, we introduce the micromorphic-type variable $\b{g}$ which is supposed to carry the gradient effects of $\varphi$.
Motivated by the seminal work of Forest~\cite{Forest2009}, we consider the associated micromorphic-type potential i.e. $\Psi_{_\text{\!MM}}=\Psi_{_\text{\!MM}}(\varphi,\b{g})$ to be a function of the generalized relative (strain-like) variable, defined by $\b{e} = \nabla \varphi - \b{g}$ and $\nabla \b{g}$, thereby introducing a coupling between original and micromorphic-type variables.
The additional contributions to the potential $\Psi_{_\text{\!MM}}$ are quadratic functions of $\b{e}$ and $\nabla \b{g}$, and the total potential is now
\begin{equation}
\begin{aligned}
&\Psi = \Psi_{0} + \Psi_{_\text{\!MM}}\,,
\end{aligned}
\label{eq:overall-energy-1}
\end{equation}
where
\begin{equation}
\begin{aligned}
&\Psi_{0} = \f{k}{2} |\nabla \varphi |^2 + F \varphi\,,\\
&\Psi_{_\text{\!MM}} = \f{1}{2} \lb \nabla \varphi - \b{g} \rb \cdot \b{H} \cdot \lb \nabla \varphi - \b{g} \rb+ \f{1}{2}  \b{g} \cdot  \b{K} \cdot  \b{g} + \f{1}{2} \nabla \b{g} :  \v{A} : \nabla \b{g}  \,,
\end{aligned}
\label{eq:overall-energy-2}
\end{equation}
with $\b{H}$ the micromorphic-type (second-order) coupling tensor, $\b{K}$ the micromorphic-type (second-order) tensor and $\v{A}$ the micromorphic-type (fourth-order) ``stiffness'' tensor.
Note that $\b{H}$, $\b{K}$ and $\v{A}$ are defined to be symmetric and positive-definite; that is, for any vector $\b{a}$ and second-order tensor $\b{A}$, $\lb\b{H}\cdot\b{a}\rb \cdot \b{a} \geq H_0 |\b{a}|^2 $, $\lb\b{K}\cdot\b{a}\rb \cdot \b{a} \geq K_0 |\b{a}|^2 $ and $\lb\v{A}:\b{A}\rb : \b{A} \geq A_0 |\b{A}|^2$ where $H_0$, $K_0$ and $A_0$ are positive constants.
The micromorphic-type coupling tensor $\b{H}$ enforces consistency between the original and micromorphic-type quantities.
In the discrete case, we will choose $\b{H}$ such that it takes the upwinding effect into account.
On the other hand, the micromorphic-type ``stiffness'' tensor $\v{A}$ defines the scale at which micromorphic effects are significant, and by appearing in gradient terms, it determines the extent to which micromorphic effects influence the solution.
One of the key features of our methodology is that, depending on the problem of interest, the influence of micromorphic effects can be completely controlled through the choice of $\b{H}$, $\b{K}$ and $\v{A}$ in order to obtain stable and accurate solutions.

\noindent
\textit{\textbf{Remark:} The choice $\v{A}=\v{O}$, $\b{K}=\bnull$ and $\b{H}=p\b{I}$ reduces the MMAD formulation to the MZAD approach described above.}

The solution variables in our problem are now $\lc \varphi , \b{g} \rc$.
Conceptually, the next step is to define the associated functional J.
This is given by
\begin{equation}
J(\varphi,\b{g}) = \int_{\cc{B}} \lb \Psi_{0}(\varphi) + \Psi_{_\text{\!MM}}(\varphi,\b{g})  \rb \d v\,,
\label{eq:new-energy}
\end{equation}
the minimizer of which corresponds to the weak problem
\begin{equation}
\begin{aligned}
&\int_{\cc{B}} k \nabla \varphi \cdot \nabla \delta \varphi \, \d v + \int_{\cc{B}}   \Big[ \b{H} \cdot  \lb \nabla \varphi - \b{g} \rb \Big] \cdot \nabla \delta \varphi \, \d v = \int_{\cc{B}} F \,  \delta \varphi \, \d v\,,\\
& \int_{\cc{B}} \Big[ -\b{H} \cdot  \lb \nabla \varphi - \b{g} \rb + \b{K} \cdot \b{g} \Big]  \cdot \delta \b{g} \, \d v  + \int_{\cc{B}} \Big[ \v{A} : \nabla \b{g} \Big] : \nabla \delta \b{g} \, \d v = 0\,.
\end{aligned}
\label{eq:weak-MMAD-0}
\end{equation}
Finally, we recover the MMAD formulation for the full convection-diffusion problem by adding the first two terms of Eq.~\eqref{eq:weak-tmp} to Eq.~\eqref{eq:weak-MMAD-0}.
Thus, we arrive at the MMAD problem in which we seek $\varphi \in \Phi$ and $\b{g} \in G$ which satisfy
\begin{equation}
\begin{aligned}
&\int_{\cc{B}}   \f{\p\varphi}{\p t} \delta \varphi \, \d v + \int_{\cc{B}} \lb \b{u} \cdot \nabla \varphi \rb \, \delta \varphi \, \d v + \int_{\cc{B}} k \nabla \varphi \cdot \nabla \delta \varphi \, \d v  + \int_{\cc{B}}   \Big[ \b{H} \cdot  \lb \nabla \varphi - \b{g} \rb \Big] \cdot \nabla \delta \varphi \, \d v = \int_{\cc{B}} F \,  \delta \varphi \, \d v\,,\\
& \int_{\cc{B}} \Big[ -\b{H} \cdot  \lb \nabla \varphi - \b{g} \rb + \b{K} \cdot \b{g} \Big]  \cdot \delta \b{g} \, \d v  + \int_{\cc{B}} \Big[ \v{A} : \nabla \b{g} \Big] : \nabla \delta \b{g} \, \d v = 0\,,
\end{aligned}
\label{eq:weak-MMAD}
\end{equation}
for all $\delta\varphi \in \Phi$ and $\delta\b{g} \in G$.
Here $\Phi=\cc{H}_0^1(\cc{B})$ as before and $G = \{\b{g}\,|\, g_i \in \cc{H}^1(\cc{B})  \}$ with the norm $||\b{g}||^2_{G} = \int_{\cc{B}} \lb |\b{g}|^2 + |\nabla \b{g}|^2 \rb \, \d v$.


\section{The discrete MMAD formulation}\label{sec:Discretization}

\noindent
We start with spatial discretization of the problem.
For convenience, since time-dependence does not play a role in our analysis, we confine our attention to the steady-state convection-diffusion problem.
As before, we define $\varphi_{h}$ and $\b{g}_{h}$ as the discrete conforming approximations of $\varphi$ and $\b{g}$, respectively, with the corresponding test functions denoted by $\delta \varphi_{h}$ and $\delta \b{g}_{h}$.
The discretized form of Eq.~\eqref{eq:weak-MMAD} reads: find $\varphi_{h} \in \Phi^{h} \subset \Phi $ and $\b{g}_{h} \in G^{h} \subset G$ such that
\begin{subequations}
\begin{align}
&\!\! \int_{\cc{B}}\!\!  \lb \b{u} \cdot \nabla \varphi_{h} \rb \, \delta \varphi_{h} \, \d v +\!\! \int_{\cc{B}}\!\! k \nabla \varphi_{h} \cdot \nabla \delta \varphi_{h} \, \d v +\!\!\! \int_{\cc{B}}   \Big[ \b{H} \cdot  \lb \nabla \varphi_{h} - \b{g}_{h} \rb \Big] \cdot \nabla \delta \varphi_{h} \, \d v = \int_{\cc{B}} \!\!F \,  \delta \varphi_{h} \, \d v\,,
\label{eq:BB-first}\\
- & \int_{\cc{B}} \Big[ \b{H} \cdot  \lb \nabla \varphi_{h} - \b{g}_{h} \rb \Big]  \cdot \delta \b{g}_{h} \, \d v  + \int_{\cc{B}} \lb \b{K} \cdot \b{g}_{h} \rb  \cdot \delta \b{g}_{h} \, \d v + \int_{\cc{B}} \Big[ \v{A} : \nabla \b{g}_{h} \Big] : \nabla \delta \b{g}_{h} \, \d v = 0\,,
\label{eq:BB-second}
\end{align}
\label{eq:weak-discretized}
\end{subequations}
for all $\delta \varphi_{h} \in \Phi^{h}$ and $\delta \b{g}_{h} \in G^{h}$.

\subsection{Well-posedness and convergence analysis}
\noindent
We define $\h{\Phi}=\Phi\times G$, with the norm $||\bullet||_{\h{\Phi}}$ defined by $||\h{\b{\varphi}}||^2_{\h{\Phi}} = ||\varphi||^2_{\Phi} + ||\b{g}||^2_{G}$ for $\h{\b{\varphi}}:=(\varphi,\b{g})\in \h{\Phi}$.
Then the fully discrete version of Eq.~\eqref{eq:weak-discretized} can be written as follows:
find $\h{\b{\varphi}}_{h} \in \h{\Phi}^{h}$such that 
\begin{equation}
B\lp\h{\b{\varphi}}_{h},\delta \h{\b{\varphi}}_{h}\rp =\ell(\delta \h{\b{\varphi}}_{h}) \,,
\label{eq:B-l-first}
\end{equation}
for all $\h{\b{\varphi}}_{h} \in \h{\Phi}^{h}$, where 
\begin{equation}
\begin{aligned}
&B\lp\h{\b{\varphi}}_{h},\delta \h{\b{\varphi}}_{h}\rp =\!\!\!\! \int_{\cc{B}} \lb  \lb \b{u}\cdot\nabla\varphi_{h}\rb \delta \varphi_{h} \!+\!  k \nabla \varphi_{h} \!\cdot\! \nabla \delta \varphi_{h} + \Big[ \b{H} \!\cdot\!  \lb \nabla \varphi_{h}\! -\! \b{g}_{h} \rb \Big] \!\cdot\! \lb \nabla \delta \varphi_{h} \!-\! \delta \b{g}_{h} \rb  \!+\! \lb\b{K}\cdot\b{g}_{h}\rb\cdot \delta \b{g}_{h} \!+\! \lb\v{A}: \nabla \b{g}_{h}\rb \!:\! \nabla \delta \b{g}_{h} \rb \, \d v\,,\\
&\ell(\delta \h{\b{\varphi}}_{h}) =  \int_{\cc{B}}  F \delta \varphi_{h}\, \d v\,.
\end{aligned}
\label{eq:B-l}
\end{equation}
To show the well-posedness of our approach, we need to show that $B(\bullet,\bullet)$ is coercive and continuous, and $\ell(\bullet)$ is continuous.
That is,
\begin{equation}
\begin{aligned}
&B(\varphi_{h}, \b{g}_{h};\varphi_{h}, \b{g}_{h}) \geq M\; ||\h{\b{\varphi}}_{h}||^2_{\h{\Phi}} \,,\\
&B(\varphi_{h}, \b{g}_{h};\delta \varphi_{h}, \delta \b{g}_{h}) \leq m\; ||\h{\b{\varphi}}_{h}||_{\h{\Phi}}\; ||\delta\h{\b{\varphi}}_{h}||_{\h{\Phi}}\,,\\
&\ell(\varphi_{h}, \b{g}_{h}) = c  \; ||\h{\b{\varphi}}_{h}||_{\h{\Phi}} \,,
\end{aligned}
\label{eq:B-l}
\end{equation}
where $M$, $m$ and $c$ are positive constants.

To show the coercivity of $B(\bullet,\bullet)$, we start with 
\begin{equation}
B(\h{\b{\varphi}}_{h}, \h{\b{\varphi}}_{h}) = \int_{\cc{B}} \lb  \lb \b{u}\cdot\nabla\varphi_{h} \rb \varphi_{h} + k |\nabla \varphi_{h}|^{2} + \lb\b{H}\cdot \lb \nabla \varphi_{h}\! -\! \b{g}_{h} \rb \rb \cdot \lb \nabla  \varphi_{h} \!-\! \b{g}_{h} \rb + \lb\b{K}\cdot\b{g}_{h}\rb\cdot\b{g}_{h}  + \lb\v{A}: \nabla \b{g}_{h}\rb \!:\! \nabla \b{g}_{h} \rb \d v\,.\\
\label{eq:B-star}
\end{equation}
Using the divergence theorem we write, for $\varphi_{h}\in \Phi^{h}$,
\begin{equation}
\int_{\cc{B}} \lb \b{u}\cdot\nabla\varphi_{h}\rb \varphi_{h}\, \d v  = 
\int_{\p\cc{B}} \varphi_{h}^2\, \b{u} \cdot \b{n} \,\d a  -
\int_{\cc{B}} \varphi_{h} \bigg[ \underbrace{\lb\nabla \cdot \b{u}\rb}_{=0} \varphi_{h} + \b{u} \cdot \nabla \varphi_{h} \bigg] \,\d v\,.
\label{eq:lhs}
\end{equation}
Now
\begin{equation}
\int_{\cc{B}} \lb \b{u}\cdot\nabla\varphi_{h}\rb \varphi_{h}\, \d v  = 
\f{1}{2} \int_{\p\cc{B}} \varphi_{h}^2\, \b{u} \cdot \b{n} \,\d a  = 0\,,
\label{eq:boundary}
\end{equation}
since $\varphi_{h}=0$ on $\p\cc{B}$.
Then from Eq.~\eqref{eq:B-star} and using Eq.~\eqref{eq:boundary}
\begin{equation}
B(\h{\b{\varphi}}_{h}, \h{\b{\varphi}}_{h}) 
\geq 
\int_{\cc{B}} 
\lb  
k |\nabla \varphi_{h}|^{2} + 
H_{0}| \nabla  \varphi_{h} \!-\! \b{g}_{h}|^2  + 
K_{0}| \b{g}_{h}|^2  + 
A_{0}| \nabla \b{g}_{h} |^{2}
\rb \d v\,.
\label{eq:B-star-star}
\end{equation}

To proceed, we utilize the inequality
\begin{equation}
2ab\leq \varepsilon a^2 + \f{1}{\varepsilon}b^2 
\quad \Longrightarrow \quad
\lb a-b \rb^2 \geq \lb 1-\varepsilon \rb a^2 + \lb 1 - \varepsilon^{-1} \rb b^2\,,
\label{eq:eps-1}
\end{equation}
which holds for any $\varepsilon$.
Then we have
\begin{equation}
\begin{aligned}
| \nabla\varphi_h - \b{g}_h |^2 
=
|\nabla \varphi_{h}|^2 - 2 \nabla \varphi_{h} \cdot \b{g}_{h} + |\b{g}_{h}|^2
&\geq 
|\nabla \varphi_{h}|^2 - 2 |\nabla \varphi_{h}|\, |\b{g}_{h}| + |\b{g}_{h}|^2\\
&\geq 
\lb 1 - \varepsilon \rb |\nabla\varphi_{h}|^2  + \lb 1-\varepsilon^{-1} \rb |\b{g}_{h}|^2 \,.
\end{aligned}
\label{eq:ineq}
\end{equation}
Consequently,  Eq.\eqref{eq:B-star-star} becomes
\begin{equation}
\begin{aligned}
B(\h{\b{\varphi}}_{h}, \h{\b{\varphi}}_{h}) 
&\geq 
\int_{\cc{B}} 
\bigg[  
k |\nabla \varphi_{h}|^{2} + 
H_{0}\lb 1-\varepsilon \rb |\nabla \varphi_{h}|^2 +
H_{0}\lb 1-\varepsilon^{-1} \rb |\b{g}_{h}|^2 +
K_{0}|\b{g}_{h}|^2 +
A_{0}| \nabla \b{g}_{h} |^{2}
\bigg] \d v\\
&=
\int_{\cc{B}} 
\bigg[    
\Big[ k + H_{0}\lb 1-\varepsilon \rb \Big] |\nabla \varphi_{h}|^{2} + 
\Big[ H_{0}\lb 1-\varepsilon^{-1} \rb + K_{0} \Big] |\b{g}_{h}|^2 +
A_{0}| \nabla \b{g}_{h} |^{2}
\bigg]   \d v\\
&\geq 
H_{0}\lb 1-\varepsilon \rb ||\varphi_{h} ||^2_{\Phi} + \Big[ H_{0}\lb 1-\varepsilon^{-1} \rb + K_{0} \Big] |\b{g}_{h}|^2 + A_{0}| \nabla \b{g}_{h} |^{2} \\
&\geq 
H_{0}\lb 1-\varepsilon \rb ||\varphi_{h} ||^2_{\Phi} + \text{min} \Big[ H_{0}\lb 1-\varepsilon^{-1} \rb + K_{0} , A_0 \Big] ||\b{g}_{h} ||^2_{G}\,,
\end{aligned}
\label{eq:B-star-star-star}
\end{equation}
choosing $\varepsilon$ such that $0<\varepsilon<1$, and $\b{K}$ such that $K_{0}+H_{0}\lb1-\varepsilon^{-1}\rb>0$.
Then
\begin{equation}
B(\h{\b{\varphi}}_{h}, \h{\b{\varphi}}_{h}) \geq M ||\,\h{\b{\varphi}}\,||_{\h{\Phi}}^{2}\,,
\end{equation}
where 
\begin{equation}
M = \text{min} \lp  H_0 \lb1-\varepsilon \rb, \, \text{min} \lp K_{0} + H_{0} \lb1-\varepsilon^{-1} \rb,\, A_0 \rp \rp\,.
\label{eq:M}
\end{equation}



For continuity we have
\begin{equation}
B(\h{\b{\varphi}}_{h}, \delta\h{\b{\varphi}}_{h}) 
= 
\int_{\cc{B}} \bigg[ 
\lb\b{u}\cdot\nabla\varphi_{h}\rb \delta\varphi_{h} + k \nabla \varphi_{h} \cdot \delta \nabla \varphi_{h} + 
\lb\b{H}\cdot \lb \nabla \varphi_{h}\! -\! \b{g}_{h} \rb \rb \cdot \lb \nabla  \delta\varphi_{h} \!-\! \delta\b{g}_{h} \rb  + 
\lb\b{K}\cdot \b{g}_{h} \rb \cdot \delta\b{g}_{h}  + 
\lb\v{A}: \nabla \b{g}_{h}\rb \!:\! \nabla\delta \b{g}_{h} 
\bigg] \d v\,,\\
\label{eq:B-star-error}
\end{equation}
so that
\begin{equation}
\begin{aligned}
|B(\h{\b{\varphi}}_{h}, \delta\h{\b{\varphi}}_{h})| 
&\leq 
\int_{\cc{B}} 
\bigg[  
u_{\text{max}}|\nabla\varphi_{h}|\, |\delta\varphi_{h}| + 
k |\nabla \varphi_{h}| \,  | \nabla \delta \varphi_{h}| + 
H_{\text{max}}  |\nabla \varphi_{h}\! -\! \b{g}_{h}|\,  |\nabla\delta  \varphi_{h} \!-\! \delta\b{g}_{h}|  + 
K_{\text{max}}  |\b{g}_{h}|\,  |\delta\b{g}_{h}|  + 
A_{\text{max}} |\nabla \b{g}_{h}| \, |\nabla \delta\b{g}_{h}| \bigg] \d v\\
&\leq
\lb 
u_{\text{max}} ||\delta\varphi_{h}||_{L^2} + 
k \,  || \nabla \delta \varphi_{h}||_{L^2} \rb ||\nabla \varphi_{h}||_{L^2} + 
2cH_{\text{max}} ||\h{\b{\varphi}}_{h}||_{\h{\Phi}} ||\delta\h{\b{\varphi}}_{h}||_{\h{\Phi}}  + 
K_{\text{max}} ||\b{g}_{h} ||_{L^2}||\delta \b{g}_{h} ||_{L^2} + 
A_{\text{max}} ||\nabla \b{g}_{h} ||_{L^2}||\nabla \delta \b{g}_{h} ||_{L^2}  \\
&\leq
\lb u_{\text{max}} + k + 2cH_{\text{max}} + K_{\text{max}} + A_{\text{max}}  \rb ||\h{\b{\varphi}}_{h}||_{\h{\Phi}} ||\delta\h{\b{\varphi}}_{h}||_{\h{\Phi}}  \,.
\end{aligned}
\label{eq:B-star-star-error}
\end{equation}
Here $u_{\text{max}} = \underset{i}{\text{max}}|u_{i}|$, $H_{\text{max}} = \underset{ij}{\text{max}}|H_{ij}|$, $K_{\text{max}} = \underset{ij}{\text{max}}|K_{ij}|$,  $A_{\text{max}} = \underset{ijkl}{\text{max}}|A_{ijkl}|$, $c=\lb\text{vol}\, \cc{B} \rb^{1/2}$ and we have used the Cauchy--Schwarz inequality and the inequality $||f||_{L^1}\leq  c||f||_{L^2}$.
Therefore $B(\bullet,\bullet)$ is continuous with constant
\begin{equation}
m = u_{\text{max}}+k  + 2cH_{\text{max}} + K_{\text{max}}  + A_{\text{max}}\,.
\label{eq:m}
\end{equation}
The proof of continuity of $\ell(\bullet)$ is straightforward and is omitted here.

\subsection{Error analysis}
\noindent
The continuity and coercivity of $B(\bullet,\bullet)$ imply the standard finite element error estimate~\cite{Reddy1998}
\begin{equation}
||\h{\b{\varphi}} - \h{\b{\varphi}}_{h} ||_{\h{\Phi}} \leq Ch \,.
\label{eq:Ch}
\end{equation}
Thus, we have convergence at a linear rate and the constant $C$ is given by $C=m/M$, where $m$ and $M$ are the continuity and coercivity constants, respectively; see for example~\cite{Brezzi1999b}.
The original semi-discrete problem is to find $\varphi_0 \in \Phi$ such that
\begin{equation}
\underbrace{\int_{\cc{B}} \bigg[ \lb  \b{u}\cdot\nabla\varphi_0\rb \delta\varphi + k \nabla \varphi_0 \cdot  \nabla \delta\varphi \bigg] \d v}_{a\lp\varphi_0, \delta \varphi \rp} 
= 
\underbrace{\int_{\cc{B}} F \, \delta \varphi \, \d v}_{\ell\lp\delta \varphi \rp} \,.\\
\label{eq:i}
\end{equation}
On the other hand, the MMAD problem is to find $\varphi\in\Phi$ and $\b{g}\in G$ such that
\begin{equation}
a\lp\varphi, \delta \varphi \rp + b\lp\varphi,\,\b{g};\,\delta\varphi \rp = \ell(\delta\varphi)\,,
\label{eq:ii}
\end{equation}
where 
\begin{equation}
b\lp\varphi,\,\b{g};\,\delta\varphi \rp = \int_{\cc{B}} \Big[\b{H}\cdot \lb \nabla \varphi\! -\! \b{g} \rb \Big] \cdot  \nabla  \delta\varphi \,  \d v\,.
\label{eq:ii-i}
\end{equation}
Thus we have an error $\lb\varphi -  \varphi_0\rb$ between the original (or actual) and MMAD solutions.
We now estimate this error.
From Eqs.~\eqref{eq:i} and~\eqref{eq:ii} we have
\begin{equation}
a\lp\varphi-\varphi_0,\, \delta \varphi \rp + b\lp \h{\b{\varphi}},\,\delta\varphi \rp = 0\,.
\label{eq:iii}
\end{equation}
Now $a\lp\bullet,\,\bullet \rp$ is coercive, i.e.
\begin{equation}
a\lp\varphi,\, \varphi \rp \geq M_0\, ||\varphi ||^2_{\Phi}\,,
\label{eq:iv}
\end{equation}
with $M_0=k$ using Eq.~\eqref{eq:boundary}.
Next, set $\delta \varphi = \varphi - \varphi_0$ in Eq.~\eqref{eq:iii}; this gives
\begin{equation}
a\lp\varphi-\varphi_0,\, \varphi-\varphi_0 \rp + b\lp \h{\b{\varphi}},\, \varphi-\varphi_0 \rp= 0\,.
\label{eq:v}
\end{equation}
Thus, from Eqs.~\eqref{eq:iv} and ~\eqref{eq:v}, the modelling error is
\begin{equation}
M_0\; ||\varphi-\varphi_0 ||^{2}_{\Phi} 
\leq 
b\lp \h{\b{\varphi}},\,\varphi_0-\varphi \rp
\leq
\beta_0 \;||\h{\b{\varphi}}||_{\h{\Phi}}\; ||\varphi - \varphi_0 ||_{{\Phi}}\,,
\label{eq:vi}
\end{equation}
where $\beta_0=H_{\text{max}}$.
Hence,
\begin{equation}
||\varphi - \varphi_0 ||_{\Phi} \leq \f{\beta_0}{M_0} ||\h{\b{\varphi}}||_{\h{\Phi}}\,.
\label{eq:vii}
\end{equation}
Assuming that $k$ is fixed or given (either an actual or artificial diffusivity), we can control the ``error'' between the actual solution and the MMAD solution through the choice of $\b{H}$.
Finally, combining Eqs.~\eqref{eq:vii} and~\eqref{eq:Ch} we have, for the error between the original solution and its MMAD finite element approximation,
\begin{equation}
||\varphi_0 - \varphi_h ||_{\Phi} \leq  ||\varphi_0 - \varphi ||_{\Phi} + ||\varphi - \varphi_h ||_{\Phi} \leq \f{\beta_0}{M_0} \,. 
\label{eq:viii}
\end{equation}

\noindent
\textit{\textbf{Remark:} 
Note that the finite element error estimate holds in the pure advection limit $k=0$.
}

\noindent
\textit{\textbf{Remark:} 
The modelling error depends on $\beta_0/M_0=H_{\text{max}}/k$. Thus it can be controlled in the advective limit through an appropriate choice for $\b{H}$. The choices for $\b{H}$, $\b{K}$ and $\v{A}$ will be discussed in the following sub-section.
}




\subsection{The choice of $\b{H}$, $\b{K}$ and $\v{A}$}\label{sec:H-A}
\noindent
The next step to complete our approach is to elaborate on the choice of $\b{H}$, $\b{K}$ and $\v{A}$.
Following the SUPG method~\cite{Brooks1982}, we set
\begin{equation}
\b{H} = \o{k}\, \h{\b{u}} \dyad \h{\b{u}}\,,
\label{eq:parameters}
\end{equation}
where $\h{\b{u}} = \b{u}/|\b{u}|$.
The parameter $\o{k}$, for example for a two-dimensional problem, is defined as
\begin{equation}
\begin{aligned}
&\o{k} = \sum_{i=1}^{\text{PD}} u_{i} h_{i} {\gamma}_{i} / 2\,,\qquad \text{with} \qquad
{\gamma}_{i} = \coth(\alpha_i)-1/\alpha_i\,,\qquad
\alpha_i=u_i h_i/2D\,,\qquad
u_i = \b{e}_i^T \cdot \b{u}\,,
\end{aligned}
\label{eq:KK}
\end{equation}
with PD being the problem dimension, and $\b{e}_{i}$ the unit vectors in the finite element natural coordinate system.
The tensors $\b{K}$ and $\v{A}$ are assumed to be of the form $\b{K}=\t{k}\b{I}$ and $\v{A}=\t{k}\v{I}$ with $\b{I}$ and $\v{I}$ being the second- and fourth-order identity tensors, respectively.
In the examples that follow we will choose either $\t{k}=1$ for problems of convection-diffusion and $\t{k}=0$  for problems of pure advection.
This selection can be compactly expressed as $\t{k}=\operatorname{sgn}(D)$.

For the choice~\eqref{eq:parameters} of $\b{H}$, we see that $\b{a}\cdot \b{H} \cdot \b{a} = \o{k} \lb \h{\b{u}} \cdot \b{a} \rb^2 \geq 0$; that is, $\b{H}$ is positive semi-definite since $\b{a}\cdot \b{H} \cdot \b{a}=0$ for $\h{\b{u}}=\bnull$ or $\h{\b{u}} \perp\b{a}$.
Likewise, the choices $\b{K}=\bnull$ and $\v{A}=\v{O}$ are trivially non-positive-definite.
These choices do not meet the \emph{sufficient} conditions for well-posedness, but as will be seen, all lead to stable and accurate numerical solutions.

\section{Numerical examples}\label{sec:result}
\noindent
In this section we evaluate the performance the MMAD formulation through a set of numerical examples.
To that end, we consider examples presented in~\cite{Brooks1981}.
The numerical schemes considered in what follows are the
\begin{itemize}
\item standard Galerkin finite element method (FEM)
\item streamline upwind Petrov--Galerkin method (SUPG)
\item mean zero artificial diffusion method (MZAD)
\item micromorphic-based artificial diffusion method (MMAD)
\end{itemize}
It was found that the streamline upwind method and the variational multi-scale method yield somewhat similar results to SUPG.
Results using these methods are therefore omitted.
To set the stage, the performance of different numerical schemes are compared against each other in a one-dimensional setting. 
Two different problems are considered; steady-state convection-diffusion and transient advection.
This is then followed by similar studies in a two-dimensional setting with different flow types.
Finally, a heat convection-conduction problem is solved using our proposed methodology.
As shown later, in order to evaluate the accuracy of each numerical scheme, two different relative error norms are reported as
\begin{equation}
\|e\|_2 = \f{\|\varphi - \varphi_h\|_2}{\| \varphi\|_2}\,, 
\qquad\qquad
\text{and}
\qquad\qquad
\|e\|_{\infty} = \f{\|\varphi - \varphi_h\|_{\infty}}{\| \varphi\|_{\infty}}\,, 
\label{eq:error-num}
\end{equation}
with $\varphi$ being the exact solution and $\varphi_h$ the discretized numerical solution.
The norms are defined as
\begin{equation}
\|\varphi\|_2 = \lp \int_{\cc{B}} | \varphi(x)|^2 \, \d v \rp ^{1/2}\,, 
\qquad\qquad
\text{and}
\qquad\qquad
\|\varphi\|_{\infty} = \text{max} \{ \varphi(x), x \in \cc{B} \}  \,.
\label{eq:error-num}
\end{equation}
These two error norms highlight the overall and point-wise errors.
Moreover, for transient advection problems, we examine the performance of the numerical schemes in maintaining the maximum value of the solution throughout a certain number of time steps.
Thus, Eq~\eqref{eq:error-num} is modified to
\begin{equation}
\|e_{\text{t}}\|_2 = \f{\lp \displaystyle\sum_{\text{n}=0}^{\#ts} |\varphi^{\text{max}}_n - \varphi_{h_{n}}^{\text{max}}|^2\rp^{1/2}}{\lp \displaystyle\sum_{\text{n}=0}^{\#ts}|\varphi^{\text{max}}_n|^2 \rp^{1/2}}  \,, 
\qquad\qquad
\text{and}
\qquad\qquad
\|e_{\text{t}}\|_{\infty} = \f{\text{max} \{  |\varphi^{\text{max}}_n - \varphi_{h_{n}}^{\text{max}}|\}}{\text{max} \{ |\varphi^{\text{max}}_n|\}}\,, 
\label{eq:error-num-2}
\end{equation}
with $\varphi^{\text{max}}_n = \text{max} \{ \varphi(x,t_n), x \in \cc{B}, t_n \in \lb0,T \rb \}$, the subscripts ``t'' referring to the relative error through time steps, ``$n$'' referring to the time step, and the maximum function taking the maximum value over all time steps.
The max function in the second error refers to the maximum over all time steps.
All numerical results are obtained using our in-house finite element code.
Continuous piecewise-linear approximations for all fields are considered for the finite element approximations through all examples.



\subsection{1D numerical examples}\label{sec:1D-numerical-example}

\subsubsection{Steady-state convection-diffusion}

\begin{figure}[b!]
\centering
\includegraphics[width=1.0\textwidth]{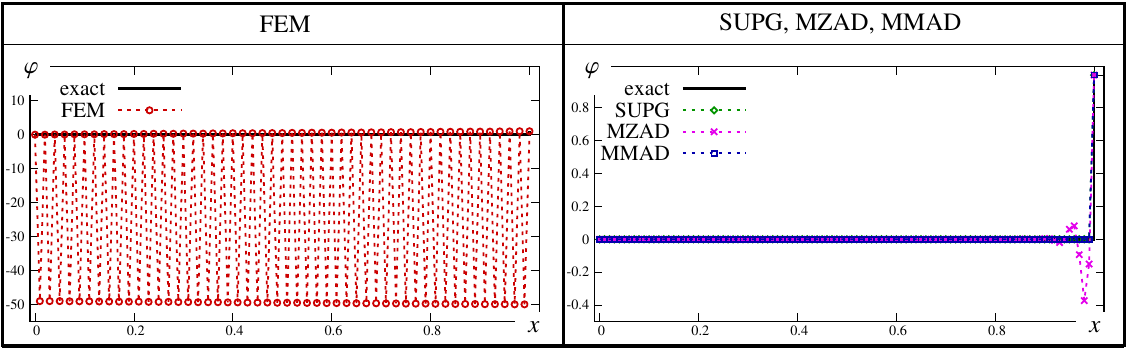}
\caption{
Comparison of various numerical methodologies for solving 1D convection-dominated problems at element Peclet number $10^{6}$.
The solid black line represents the exact solution and the dashed lines (with points on top) represent the solutions obtained from different numerical schemes.
}
\label{fig:1D-Convection-Diffusion}
\end{figure}

\begin{table}[b!]
\centering
\begin{tabular}{ |p{2.8cm}|p{2.8cm}|p{2.8cm}|p{2.8cm}|p{2.8cm}|}\hline
               & FEM         & SUPG          & MZAD       & MMAD\\\hline
$\|e\|_2     $ & $350.07$    & $\approx 0$ & $0.4265$ & $\approx 0$ \\ \hline
$\|e\|_\infty$ & $49.99$     & $\approx 0$ & $0.3731$ & $\approx 0$ \\ \hline
\end{tabular}
\caption{ 
Error of each numerical scheme for solving the 1D convection-dominated problem at element Peclet number $10^{6}$.
The $L^2$- and max-norms of the relative error associated with each numerical method are shown.
}
\label{tab:conv-diff-1D}
\end{table}

\noindent
In this example, a steady-state convection-diffusion problem is solved in a 1D domain of unit length.
The transient term $\p \varphi /\p t$ and the source term $F$ in Eq.~\eqref{eq:strong-1} are zero. 
The domain is subject to the boundary conditions $\varphi=0$ at $x=0$ and $\varphi=1$ at $x=1$.
The domain is discretized into $100$ linear finite elements and the problem is convection-dominated with element Peclet number $\text{Pe}_h=10^{6}$.
Fig.~\ref{fig:1D-Convection-Diffusion} shows the exact solution versus the solution obtained from various numerical schemes.
The solid black line represents the exact analytical solution whereas the dashed lines (with points on top) represent the numerical solutions. 
As expected, we observe that the standard Galerkin FEM fails to provide a stable solution.
The solution using the MZAD method exhibits some oscillations close to the boundary $x=1$.  
It is observed that the SUPG and MMAD solutions do not suffer from any instability and are capable of capturing the sharp gradient in the solution.
Table~\ref{tab:conv-diff-1D} shows the $L^2$-norm and max-norm of the relative error associated with each numerical scheme.
The reason for presenting both $L^2$ and max norms is to assess whether smoother solutions always lead to higher accuracy; that is, whether some oscillations leading to large max-norms might lead to better accuracy with lower $L^2$-norm.
In this example, clearly both SUPG and MMAD methods provide the most accurate solutions.


\subsubsection{Transient advection}

\begin{figure}[h!]
\centering
\includegraphics[width=1.0\textwidth]{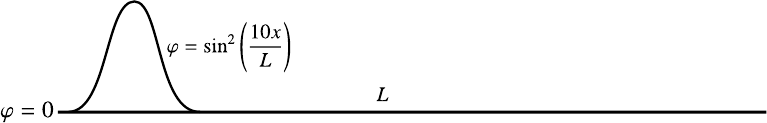}
\includegraphics[width=1.0\textwidth]{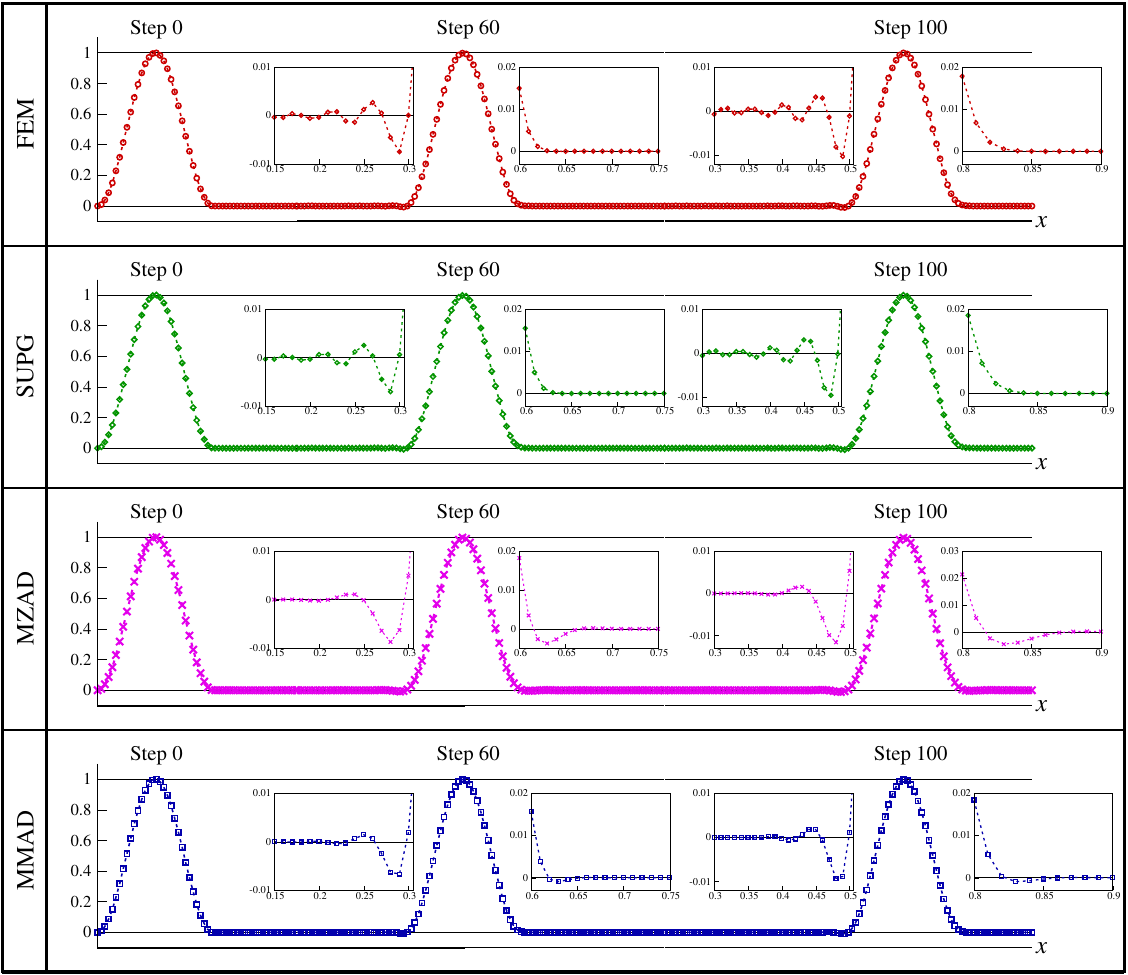}
\caption{ 
Problem definition together with initial and boundary conditions for a 1D transient advection problem (top).
Comparison of various methodologies for solving transient advection problems (bottom).
Three snapshots of the solution are provided for each case; solution at step 0, step 60 and step 100.
The zoom boxes for the intermediate steps are provided which magnify the solution at the left and right part of the sine hill.
}
\label{fig:1D-Advection-Transient-Transport}
\end{figure}

\begin{table}[b!]
\centering
\begin{tabular}{ |p{1.5cm}|p{1.5cm}|p{2.60cm}|p{2.60cm}|p{2.60cm}|p{2.60cm}| }\hline
step &                                 & FEM       & SUPG          & MZAD       & MMAD\\\hline
\multirow{2}{*}{60}   & $\|e\|_2     $ & $0.00529$  & $0.00527$     & $0.00525$   & $0.00521$ \\ 
                      & $\|e\|_\infty$ & $0.00995$  & $0.00954$     & $0.00998$   & $0.00914$ \\ \hline
\multirow{2}{*}{100}  & $\|e\|_2     $ & $0.00813$  & $0.00805$     & $0.00799$   & $0.00782$ \\ 
                      & $\|e\|_\infty$ & $0.01203$  & $0.01161$     & $0.01205$   & $0.01148$ \\ \hline
\end{tabular}
\caption{ 
Error of each numerical schemes for solving the 1D transient advection problem at time steps $60$ and $100$.
The $L^2$- and max-norms of the relative error associated with each numerical method is shown, for steps 60 and 100.
}
\label{tab:tran-adv-1D}
\end{table}

\begin{table}[b!]
\centering
\begin{tabular}{ |p{2.8cm}|p{2.8cm}|p{2.8cm}|p{2.8cm}|p{2.8cm}|}\hline
               & FEM         & SUPG          & MZAD       & MMAD\\\hline
$\|e_{\text{t}}\|_2     $ & $0.00136$     & $0.00260$ & $0.00241$ & $0.00111$ \\ \hline
$\|e_{\text{t}}\|_\infty$ & $0.00338$     & $0.00538$ & $0.00503$ & $0.00235$ \\ \hline
\end{tabular}
\caption{ 
Comparison of different numerical schemes in maintaining the magnitude of the sine hill.
Point-wise normalized errors associated with each scheme is calculated during $120$ time steps.
The $L^2$- and max-norms of the relative error associated with each numerical method is shown.
}
\label{tab:tran-adv-top-1D}
\end{table}

\noindent
Here we consider a transient advection problem.
The source term $F$ in Eq.~\eqref{eq:strong-1} is zero and the diffusion coefficient $D=10^{-6}$ leading to a convection-dominated problem.
In this example a sine hill function $\sin^{2}\lp10x/L\rp$ is transported over the domain, as shown in Fig.~\ref{fig:1D-Advection-Transient-Transport} (top).
Three snapshots of the solution associated with each numerical scheme are provided;
step $0$ showing the initial state, step $60$ and step $100$.
Zoom boxes magnify the solution at the left and right part of the sine hill to provide more resolution.
The domain is discretized into $100$ linear finite elements and the Crank-Nicolson scheme is employed for time integration.
Moreover, in each graph, a solid straight line at $\varphi=1$ is shown as a reference to facilitate observation of over-smoothing or over-damping of the sine hill.

In the first row, the standard FEM is observed to perform well in terms of maintaining the hill size.
However, small oscillations in the form of overshooting and undershooting are observed at the left side of the hill on each step. 
The solution by the SUPG method does not suffer from over-damping and is more or less similar to the FEM solution where small oscillations occur at the left side of the hill.
The MZAD method provides a solution that is less oscillatory on the left side of the hill.
However, some undershoots are observed on the right side of the hill.
The MMAD method also provides an accurate solution in terms of no over-damping.
Small oscillations on the left hand side and very tiny undershoots on the right hand side of the hill occur on each time step.
Such oscillations seem to disappear earlier than with the SUPG and MZAD methods.

In order to draw a more systematic comparison between the accuracy of various numerical schemes in this example, two different analyses are conducted.
In the first analysis, the solution of each numerical scheme is compared with the exact solution at time steps $60$ and $100$ and the associated $L^2$- and max-norms of the relative error are reported in Table~\ref{tab:tran-adv-1D}.
Clearly, the SUPG, MZAD and MMAD methods yield more accurate solutions than FEM with the MMAD method providing the most accurate solution.
MZAD shows a smaller $L^2$-norm error compared to SUPG while the opposite holds for the max-norm.
In the second analysis, the value of the hill top at $120$ steps is recorded and compared against the exact value $\varphi=1$.

Table~\ref{tab:tran-adv-top-1D} shows the $L^2$ and max relative error norms associated with each numerical scheme in maintaining the sine hill throughout all the time steps.
It is observed that among all the methodologies, MMAD leads to the least decay in time.
Surprisingly, the SUPG and MZAD methods perform even worse than FEM in maintaining the sine hill.


\subsection{2D numerical examples}\label{sec:2D-numerical-example}

\subsubsection{Steady-state convection-diffusion}

\paragraph{Irrotational flow}

\begin{figure}[b!]
\centering
\includegraphics[width=0.8\textwidth]{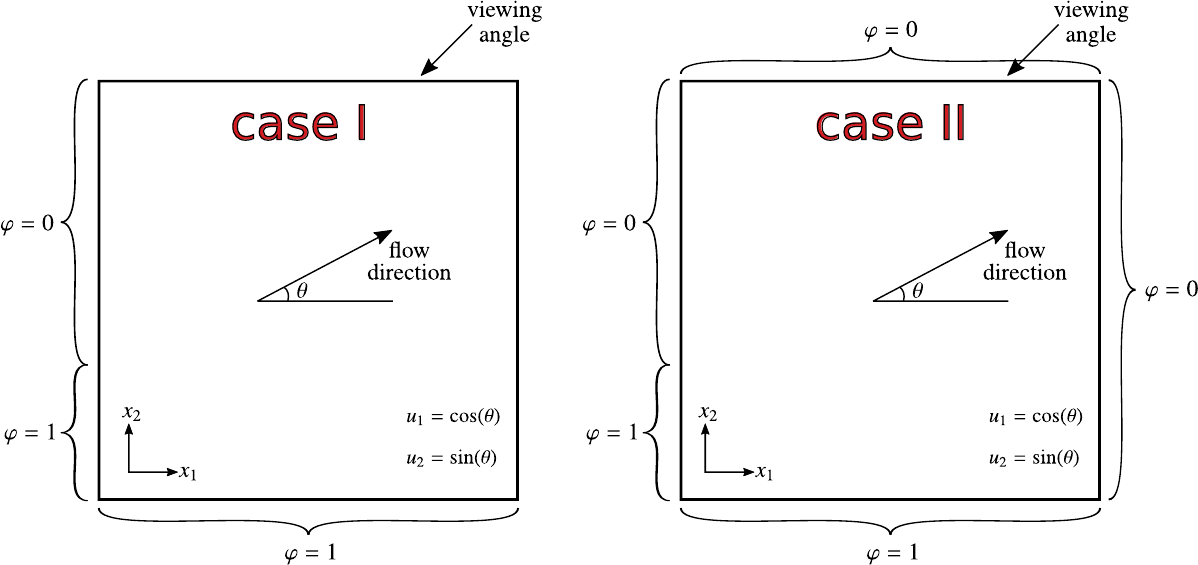}
\caption{ 
Problem definition together with boundary conditions for a 2D convection-diffusion problem.
For case I, two edges of the domain are subject to Dirichlet boundary conditions and the other edges are subject to Neumann boundary conditions.
For case II, all four edges are subject to Dirichlet boundary conditions.
}
\label{fig:2D-Convection-Diffusion-Example}
\end{figure}

\begin{figure}[h!]
\centering
\includegraphics[width=1.0\textwidth]{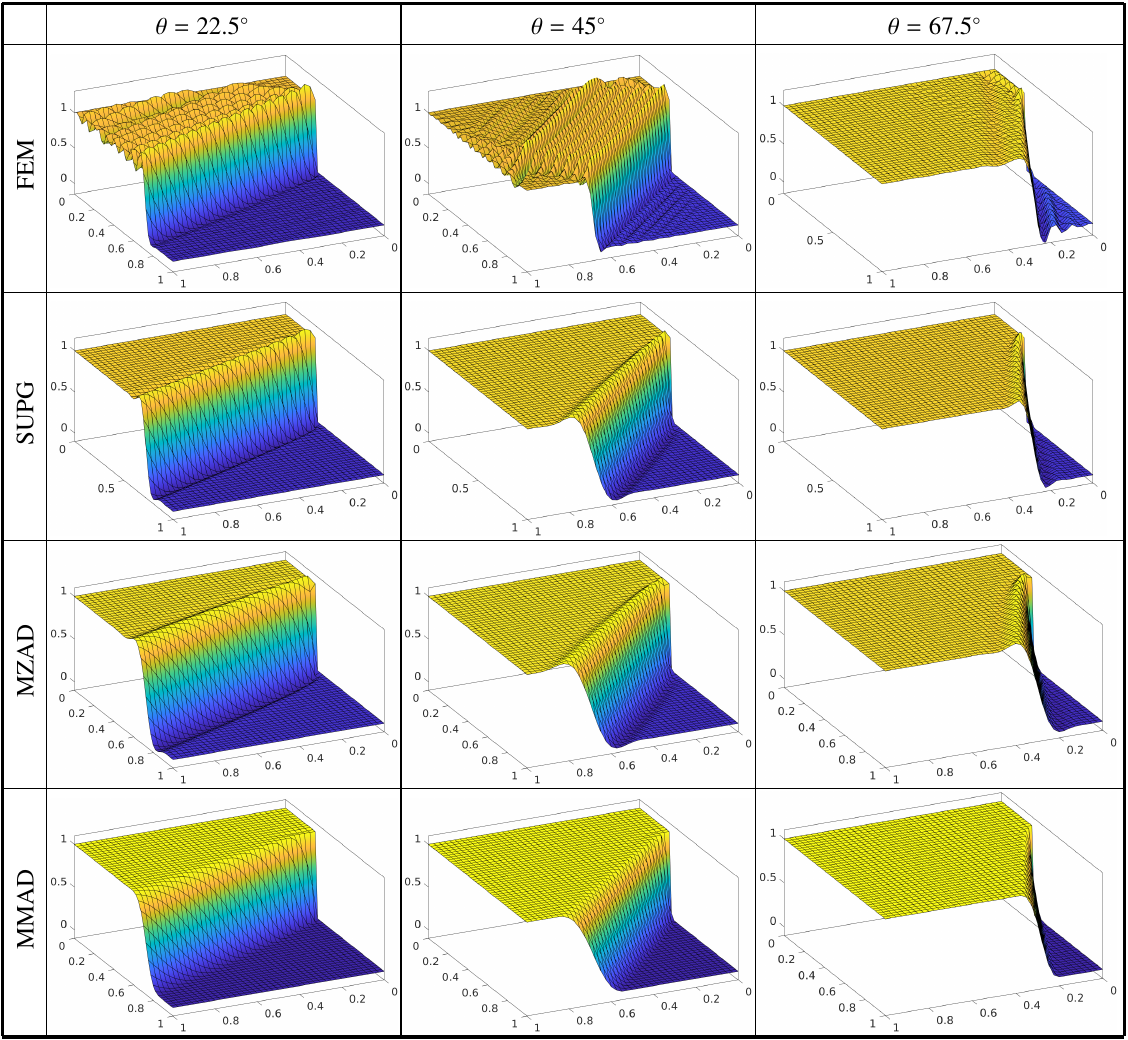}
\caption{ 
Comparison of various formulations for solving 2D convection-dominated problems at element Peclet number $10^{6}$.
The solutions are associated with test case I in Fig.~\ref{fig:2D-Convection-Diffusion-Example}.
}
\label{fig:2D-Convection-Diffusion-Free}
\end{figure}

\begin{table}[b!]
\centering
\begin{tabular}{ |p{1.5cm}|p{1.5cm}|p{2.60cm}|p{2.60cm}|p{2.60cm}|p{2.60cm}| }\hline
angle                            &                & FEM       & SUPG         & MZAD       & MMAD\\\hline
\multirow{2}{*}{$\theta=22.5$}   & $\|e\|_2     $ & $0.1179$  & $0.1013$     & $0.1107$   & $0.0947$ \\ 
                                 & $\|e\|_\infty$ & $0.6322$  & $0.5556$     & $0.5123$   & $0.4890$ \\ \hline
\multirow{2}{*}{$\theta=45$}     & $\|e\|_2     $ & $0.1041$  & $0.0617$     & $0.0628$   & $0.0573$ \\ 
                                 & $\|e\|_\infty$ & $0.4121$  & $0.3763$     & $0.3377$   & $0.3233$ \\ \hline
\multirow{2}{*}{$\theta=67.5$}   & $\|e\|_2     $ & $0.0880$  & $0.0860$     & $0.0870$   & $0.0826$ \\ 
                                 & $\|e\|_\infty$ & $0.6893$  & $0.6140$     & $0.5615$   & $0.4917$ \\ \hline
\end{tabular}
\caption{ 
Error of each numerical scheme for solving 2D convection-dominated problems at element Peclet number $10^{6}$ for test case I.
The $L^2$- and max-norms of the relative error associated with each numerical method is shown for each flow angle.
}
\label{tab:conv-diff-I-2D}
\end{table}

\begin{figure}[h!]
\centering
\includegraphics[width=1.0\textwidth]{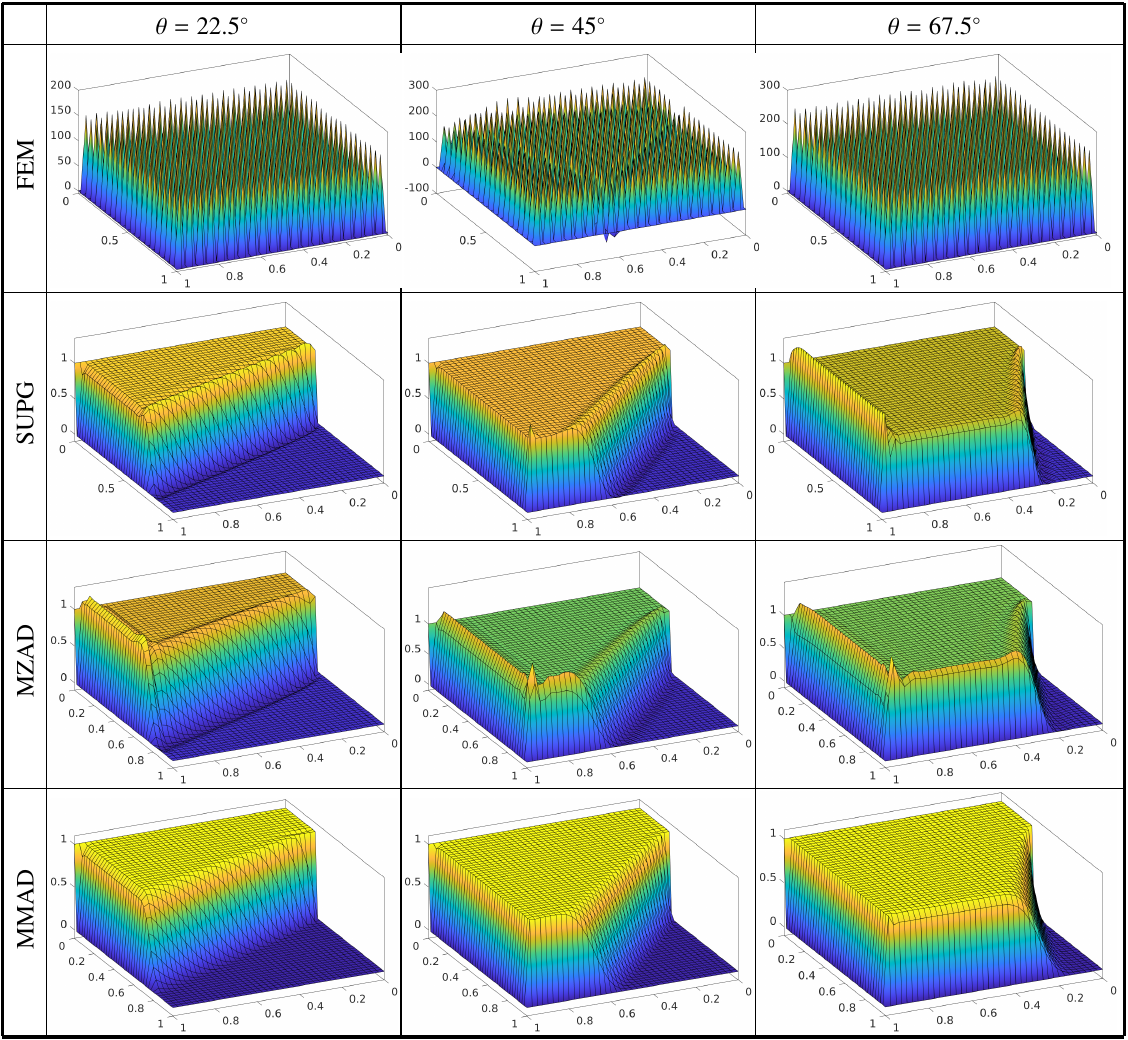}
\caption{ 
Comparison of various formulations for solving 2D convection-dominated problems at element Peclet number $10^{6}$.
The solutions are associated with test case II in Fig.~\ref{fig:2D-Convection-Diffusion-Example}.
}
\label{fig:2D-Convection-Diffusion-Fixed}
\end{figure}

\begin{table}[b!]
\centering
\begin{tabular}{ |p{1.5cm}|p{1.5cm}|p{2.60cm}|p{2.60cm}|p{2.60cm}|p{2.60cm}| }\hline
angle                            &                & FEM       & SUPG         & MZAD       & MMAD\\\hline
\multirow{2}{*}{$\theta=22.5$}   & $\|e\|_2     $ & $106.58$  & $0.1009$     & $0.1187$   & $0.0930$ \\ 
                                 & $\|e\|_\infty$ & $160.13$  & $0.5556$     & $0.6123$   & $0.4890$ \\ \hline
\multirow{2}{*}{$\theta=45$}     & $\|e\|_2     $ & $122.38$  & $0.0607$     & $0.0812$   & $0.0540$ \\ 
                                 & $\|e\|_\infty$ & $232.87$  & $0.3727$     & $0.5600$   & $0.2921$ \\ \hline
\multirow{2}{*}{$\theta=67.5$}   & $\|e\|_2     $ & $124.60$  & $0.1017$     & $0.1244$   & $0.0858$ \\ 
                                 & $\|e\|_\infty$ & $247.95$  & $0.6140$     & $0.8115$   & $0.4717$ \\ \hline
\end{tabular}
\caption{ 
Error of each numerical scheme for solving 2D convection-dominated problems at element Peclet number $10^{6}$ for test case II.
The $L^2$- and max-norms of the relative error associated with each numerical method is shown for each flow angle.
}
\label{tab:conv-diff-II-2D}
\end{table}

\noindent
In this 2D example we solve a steady-state convection-diffusion problem with irrotational flow ($\curl\, \b{u}=\bnull$). 
The transient term $\p \varphi /\p t$ and the source term $F$ in Eq.~\eqref{eq:strong-1} are zero.
The domain is a unit square $[0,1]\times[0,1]$, and is discretized into $1600$ bilinear finite elements.
We consider two different cases for the boundary conditions as shown in Fig.~\ref{fig:2D-Convection-Diffusion-Example}; case I where two edges of the domain are subject to Dirichlet boundary conditions and the other edge are subject to Neumann boundary conditions and case II where all the four edges are subject to Dirichlet boundary conditions.
The problem is convection-dominated with element Peclet number $\text{Pe}_h=10^{6}$.
The velocity components are defined as $u_{1}=\cos\theta$ and $u_{2}=\sin\theta$ with $\theta$ being the flow direction.
We consider three different flow directions of $22.5^{\circ}$, $45^{\circ}$ and $67.5^{\circ}$ for this study.

Fig.~\ref{fig:2D-Convection-Diffusion-Free} shows the solution obtained from different numerical schemes and different flow directions for case I.
Each row is associated with a numerical scheme and each column is associated with a flow direction.
As expected, the FEM solution exhibits substantial oscillatory behavior.
The SUPG scheme provides a smooth and accurate solution in the regions that are away from the shock.
Small overshoots and undershoots are observed in the vicinity of the shock.
The MZAD solution looks similar to SUPG but with slightly smoother transitions adjacent to the shock. 
Finally, the MMAD method provides a stable and oscillation-free solution for all cases.
In contrast to SUPG and MZAD, no overshoots and undershoots occur in the vicinity of the shock.

It is arguable whether the sharp but locally oscillatory SUPG and MZAD solution is more accurate or less accurate than slightly smooth MMAD solution.
To shed light on this issue, Table~\ref{tab:conv-diff-I-2D} shows a comparison of the errors associated with each numerical method and for each flow angle.
Clearly, the highest error is associated with the FEM solution.
The relatively smaller max-norm error of the MZAD method compared to the SUPG method reflects its smoother solution with lower amplitudes (but longer wave lengths) of overshoots and undershoots.
On the other hand, the smaller $L^2$-norm of the SUPG method signifies its greater overall accuracy than the MZAD method.
The MMAD method proves to be the most accurate scheme with the smallest $L^2$ and max error norms.

Next, Fig.~\ref{fig:2D-Convection-Diffusion-Fixed} illustrates the solutions of each numerical scheme for the test example II.
The FEM solution is completely unstable for all flow angles.
For the SUPG and MZAD schemes, the same overshooting and undershooting occur close to the shock for all flow directions.
In addition, for SUPG, some instabilities also occur close to the boundaries where the flow direction is $45^{\circ}$ and $67.5^{\circ}$.
Such instabilities are more significant for MZAD and they occur for all flow directions.
The solution obtained by the MMAD method proves to be stable.
No wiggles occur either in the vicinity of the shock or close to the boundary.
Table~\ref{tab:conv-diff-II-2D} is the counterpart of Table~\ref{tab:conv-diff-I-2D} for the test case II.
For this case, the SUPG method outperforms the MZAD providing both smaller $L^2$ and max error norms.
Similar to the previous example, the MMAD proves to be the most accurate scheme.


\subsubsection{Transient advection}

\paragraph{Irrotational flow}

\begin{figure}[b!]
\centering
\includegraphics[scale=0.7]{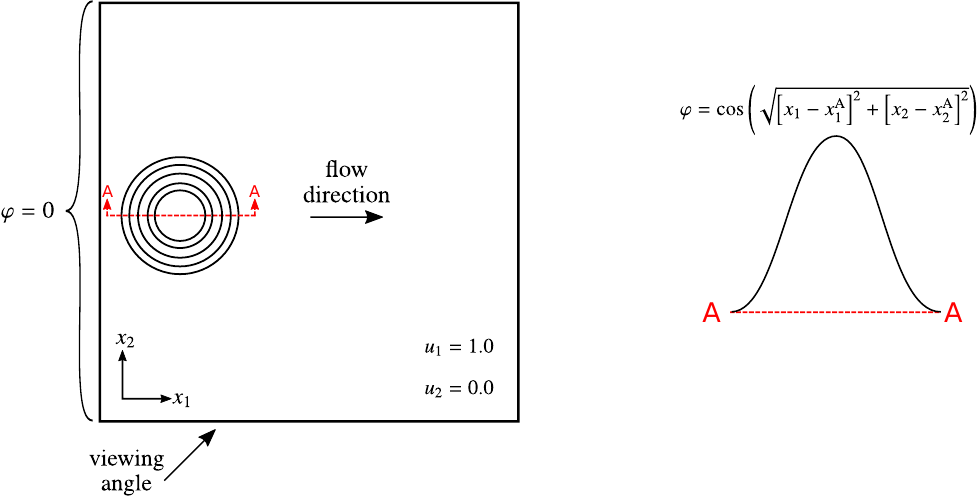}
\caption{
Problem definition together with initial and boundary conditions for a 2D transient advection problem with an irrotational flow.
}
\label{fig:2D-Advection-Transient-Example}
\end{figure}

\begin{figure}[b!]
\centering
\includegraphics[width=1.0\textwidth]{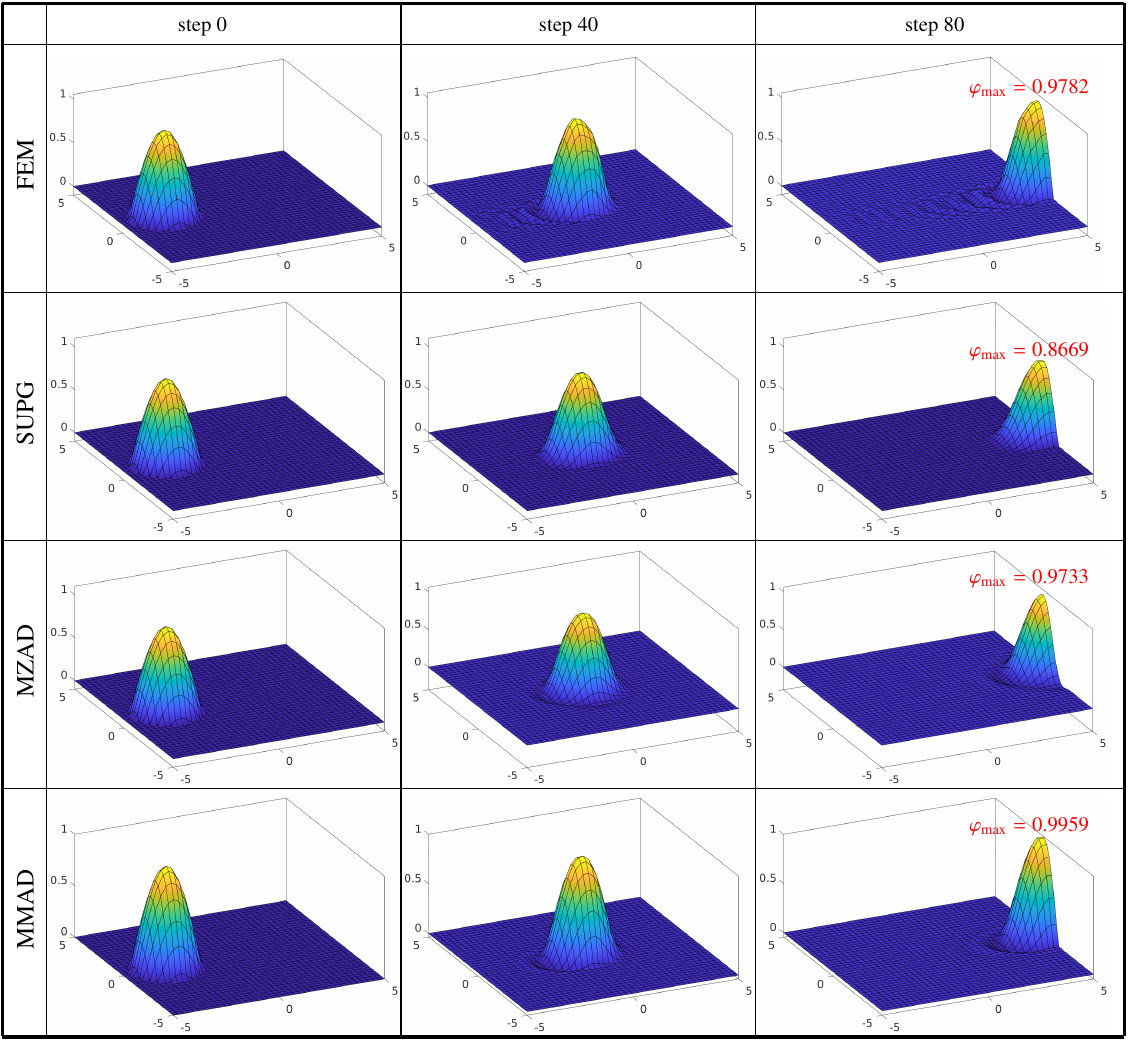}
\caption{
Comparison of various formulations for solving 2D transient advection problems with an irrotational flow.
The solutions are associated with the example described in Fig.~\ref{fig:2D-Advection-Transient-Example}.
}
\label{fig:2D-Advection-Transient-Transport}
\end{figure}

\begin{table}[b!]
\centering
\begin{tabular}{ |p{1.5cm}|p{1.5cm}|p{2.60cm}|p{2.60cm}|p{2.60cm}|p{2.60cm}| }\hline
step &                                 & FEM       & SUPG          & MZAD       & MMAD\\\hline
\multirow{2}{*}{40}   & $\|e\|_2     $ & $0.0736$  & $0.1236$     & $0.0719$   & $0.0689$ \\ 
                      & $\|e\|_\infty$ & $0.0697$  & $0.1350$     & $0.0925$   & $0.0786$ \\ \hline
\multirow{2}{*}{80}   & $\|e\|_2     $ & $0.1319$  & $0.2032$     & $0.1308$   & $0.1220$ \\ 
                      & $\|e\|_\infty$ & $0.1296$  & $0.1744$     & $0.1319$   & $0.1196$ \\ \hline
\end{tabular}
\caption{ 
Error of each numerical schemes for solving the 2D transient advection problem with irrotational flow at time steps $40$ and $80$.
The $L^2$- and max-norms of the relative error associated with each numerical method is shown.
}
\label{tab:tran-adv-irr-2D}
\end{table}

\begin{table}[b!]
\centering
\begin{tabular}{ |p{2.8cm}|p{2.8cm}|p{2.8cm}|p{2.8cm}|p{2.8cm}|}\hline
               & FEM         & SUPG          & MZAD       & MMAD\\\hline
$\|e_{\text{t}}\|_2     $ & $0.0150$     & $0.0795$ & $0.0061$ & $0.0050$ \\ \hline
$\|e_{\text{t}}\|_\infty$ & $0.0303$     & $0.1377$ & $0.0120$ & $0.0112$ \\ \hline
\end{tabular}
\caption{ 
Comparison of different numerical schemes in maintaining the magnitude of the cosine hill.
Point-wise normalized errors associated with each scheme is calculated during $80$ time steps.
The $L^2$- and max-norms of the relative error associated with each numerical method is shown.
}
\label{tab:tran-adv-irr-top-2D}
\end{table}

\noindent
In this example we consider a 2D transient advection problem with irrotational flow.
As shown in Fig.~\ref{fig:2D-Advection-Transient-Example}, a cosine hill along the rotating line AA located at the left hand side of the domain is transported towards the right hand side of the domain. 
The domain is the square $[-5,5]\times[-5,5]$.
The equation of the hill is shown in Fig.~\ref{fig:2D-Advection-Transient-Example} (right) with $x_1^{A} = -3$ and $x_2^{A} = 0$.
Dirichlet boundary conditions $\varphi=0$ are imposed on the left edge of the domain and the other three edges are subject to Neumann boundary conditions.
The velocity components are $u_{1}=1$ and $u_{2}=0$.

Fig.~\ref{fig:2D-Advection-Transient-Transport} shows the solution obtained by different numerical methods.
Three different time steps are considered corresponding to each column; time steps $0$ associated with the initial condition, step $40$ as an intermediate step and step $80$ which is intentionally chosen since it is the step at which the top of the cosine hill reaches the end of the domain.
For all numerical methods, the maximum value of the cosine hill at step $80$ is shown.
It is observed that for the FEM solution, the hill leaves an oscillatory trail behind as it transports through the domain.
Moreover, the hill height is subject to a slight decay during the hill's transport and reaches $97\%$ of its original height.
The SUPG method renders a stable smooth solution with almost no trail left behind.
However, the hill is reduced to $86.69\%$ of its original height as it reaches the right boundary of the domain.
The MZAD method also yields a stable smooth solution with a small undershoot behind the hill but maintains $97.33\%$ of the hill height during the transport.
This value for the MMAD method is $99.59\%$ which signifies the excellent accuracy of the method.

Similar to the 1D transient advection example, Table~\ref{tab:tran-adv-irr-2D} compares the accuracy of the numerical schemes at times steps $40$ and $80$.
Of particular interest is the poor performance of SUPG compared to FEM due to over smoothening.
The MZAD method performs slightly better than FEM and the MMAD method proves to be the most accurate scheme.
In the second analysis, the value of the hill top during $80$ steps is recorded and compared against the exact value $\varphi=1$.
Table~\ref{tab:tran-adv-irr-top-2D} shows the $L^2$-norm and max-norm of the error associated with each numerical scheme in maintaining the sine hill throughout all the time steps.
It is observed that among all the methodologies, MMAD is the most precise methodology rendering the least amount of decay in time.
Surprisingly, the SUPG method perform even worse than FEM in maintaining the sine hill with MZAD outperforming FEM only marginally.


\begin{figure}[b!]
\centering
\includegraphics[scale=0.7]{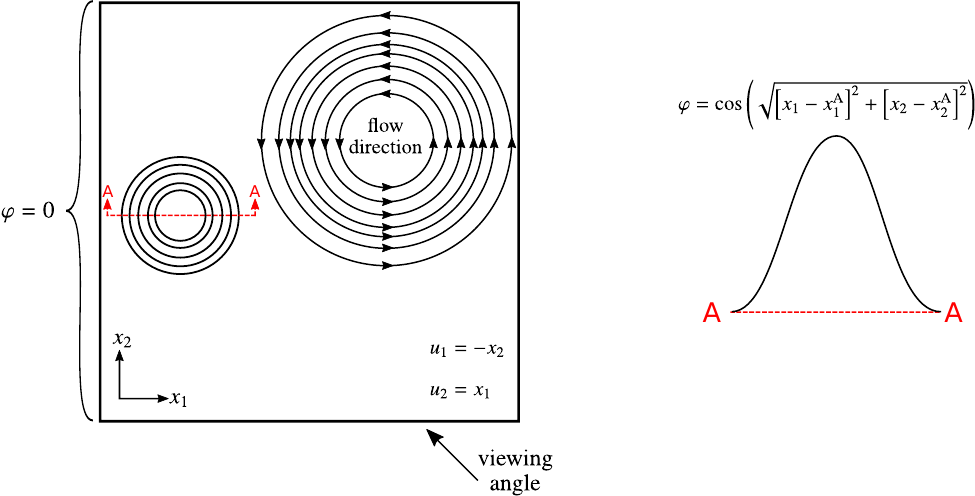}
\caption{ 
Problem definition together with initial and boundary conditions for a 2D transient advection problem with rotational flow.
The left edge is subject to Dirichlet boundary conditions whereas the other three edges are subject to Neumann boundary conditions.
A cosine hill is rotated in a square domain whose left edge is subject to Dirichlet boundary condition $\varphi=0$ and the other edges are free.
}
\label{fig:2D-Advection-Transient-Rotational-Example}
\end{figure}

\begin{figure}[b!]
\centering
\includegraphics[width=1.0\textwidth]{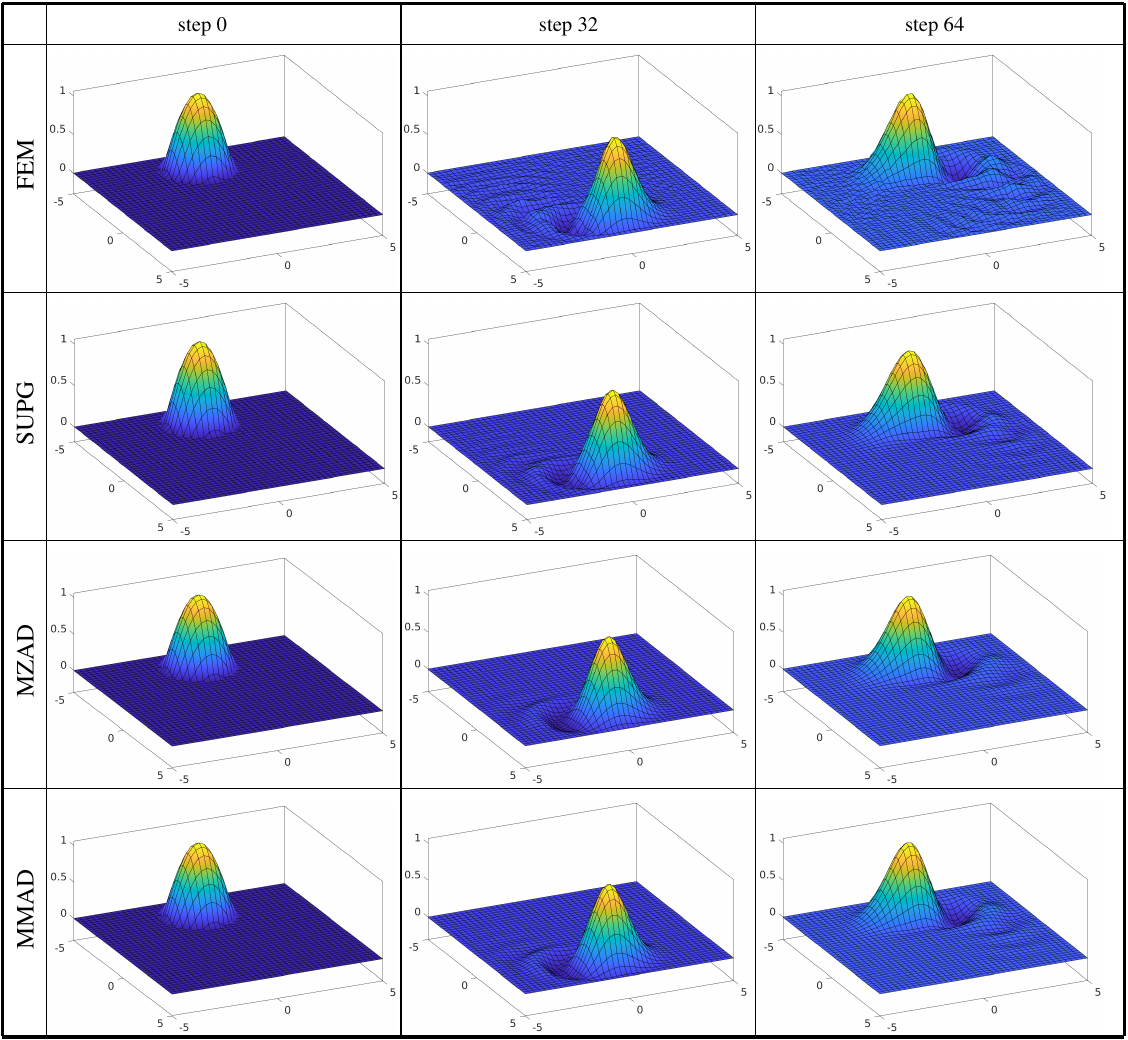}
\caption{ 
Comparison of various numerical methodologies for solving 2D transient advection problems with a rotational flow after one full rotation of the hill.
The solutions are associated with the example described in Fig.~\ref{fig:2D-Advection-Transient-Rotational-Example}.
}
\label{fig:2D-Advection-Transient-Rotational}
\end{figure}

\begin{table}[b!]
\centering
\begin{tabular}{ |p{1.5cm}|p{1.5cm}|p{2.60cm}|p{2.60cm}|p{2.60cm}|p{2.60cm}| }\hline
step &                                 & FEM       & SUPG          & MZAD       & MMAD\\\hline
\multirow{2}{*}{32}   & $\|e\|_2     $ & $0.2720$  & $0.1773$     & $0.1761$   & $0.1729$ \\ 
                      & $\|e\|_\infty$ & $0.2510$  & $0.2005$     & $0.1821$   & $0.1719$ \\ \hline
\multirow{2}{*}{64}   & $\|e\|_2     $ & $0.2865$  & $0.2609$     & $0.2611$   & $0.2564$ \\ 
                      & $\|e\|_\infty$ & $0.2640$  & $0.2818$     & $0.2530$   & $0.2463$ \\ \hline
\end{tabular}
\caption{ 
Error of each numerical schemes for solving the 2D transient advection problem with irrotational flow at time steps $32$ and $64$.
The $L^2$- and max-norms of the relative error associated with each numerical method is shown.
}
\label{tab:tran-adv-rot-2D}
\end{table}

\begin{table}[b!]
\centering
\begin{tabular}{ |p{2.8cm}|p{2.8cm}|p{2.8cm}|p{2.8cm}|p{2.8cm}|}\hline
               & FEM         & SUPG          & MZAD       & MMAD\\\hline
$\|e_{\text{t}}\|_2     $ & $0.0347$     & $0.0551$ & $0.0276$ & $0.0233$ \\ \hline
$\|e_{\text{t}}\|_\infty$ & $0.0625$     & $0.1151$ & $0.0592$ & $0.0393$ \\ \hline
\end{tabular}
\caption{ 
Comparison of different numerical schemes in maintaining the magnitude of the sine hill.
Point-wise normalized errors associated with each scheme is calculated during $64$ time steps.
The $L^2$- and max-norms of the relative error associated with each numerical method is shown.
}
\label{tab:tran-adv-rot-top-2D}
\end{table}

\paragraph{Rotational flow}
\noindent
This example involves a transient advection problem with rotational flow ($\curl\,\b{u}\neq\bnull$).
Again, the source term $F$ in Eq.~\eqref{eq:strong-1} is zero and the diffusion coefficient $D=10^{-6}$ leading to a convection-dominated problem.
As shown in Fig.~\ref{fig:2D-Advection-Transient-Rotational-Example}, a cosine hill is rotated in a square domain whose left edge is subject to Dirichlet boundary condition $\varphi=0$ and the other edges are free.
The flow velocity components are determined as $u_{1}=-x_{2}$ and $u_{2}=x_{1}$.
Fig.~\ref{fig:2D-Advection-Transient-Rotational} shows the solution of different numerical schemes after half and full rotation of the hill, corresponding to time steps $32$ and $64$, respectively.
It is observed that the FEM solution involves a very long trail on the trajectory of the hill.
SUPG provides a nice stable solution with a shorter trail compared to FEM.
The solutions by the MZAD and MMAD methods are also stable with shorter trails compared to FEM.
In order to better compare the performance of the numerical schemes for this example, Tables~\ref{tab:tran-adv-rot-2D} and~\ref{tab:tran-adv-rot-top-2D} provide a quantitative comparison between the FEM, SUPG, MZAD and MMAD methods.
In Table~\ref{tab:tran-adv-rot-2D}, the accuracy of the numerical schemes at times steps $32$ and $64$ is presented.
In contrast to the previous example with irrotational flow, the solution due to the SUPG method is more accurate than the FEM.
The MZAD solution performs slightly better than SUPG and the MMAD provides the most accurate solution.
Finally, in the last analysis, the hill decay associated with this problem is investigated for all methods.
In doing so, the maximum value of $\varphi$ associated with the top of the hill is recorded at each step during the whole rotation and compared against the exact value $\varphi=1$ throughout all the time steps.
Table~\ref{tab:tran-adv-rot-top-2D} shows the $L^2$ and max error norms associated with each numerical scheme in maintaining the cosine hill.
It is observed that among all the methodologies, MMAD is the most precise methodology rendering the least amount of decay with respect to time.
The solution due to the SUPG method appears to be over damped compared to solutions by FEM and MZAD.


\subsection{Transient heat convection-conduction}\label{sec:heat}
\begin{figure}[h!]
\centering
\includegraphics[width=1.0\textwidth]{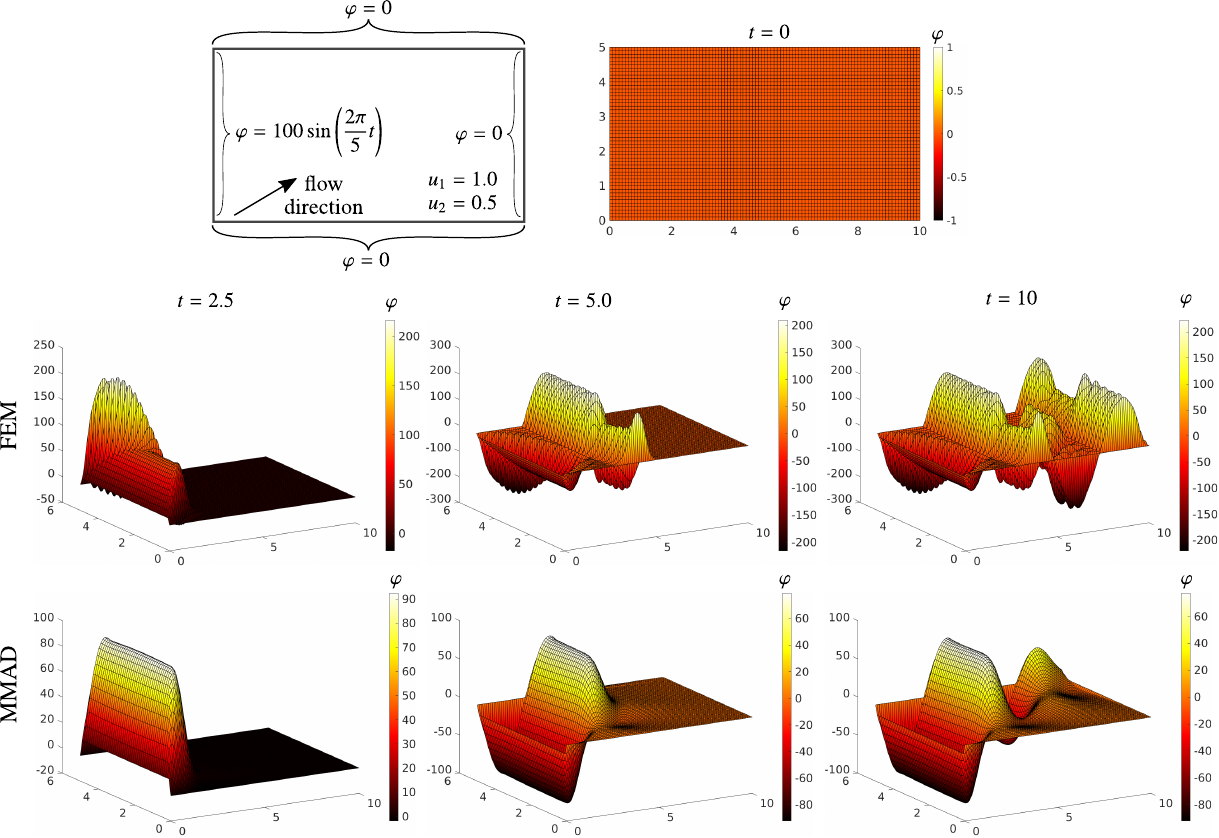}
\caption{ 
Problem definition and the solution of heat convection-conduction in a 2D rectangular domain.
The solution is obtained via the MMAD method.
The flow velocity components are determined as $u_{1}=1.0$ and $u_{2}=0.5$ and the conduction coefficient is $k=0.05$.
}
\label{fig:Heat}
\end{figure}

\noindent
As the last example, we analyze transient heat convection-conduction in a rectangular domain with time dependent boundary conditions.
That is, only the source term $F$ in Eq.~\eqref{eq:strong-1} is zero.
The variable $\varphi$ represents the temperature in this example.
As shown in Fig.~\ref{fig:Heat} top left, a rectangular domain is subject to constant temperature $\varphi=0$ on three of its edges and a time dependent sinusoidal temperature $\varphi=100\sin(2\pi/5)t$ on its left edge.
The flow velocity components are determined as $u_{1}=1.0$ and $u_{2}=0.5$ and the conduction coefficient is $D=10^{-6}$ leading to a convection-dominated problem.
For the sake of brevity, only the solutions obtained via FEM and MMAD method are shown here.
Three different snapshots of the solution at $t=2.5$, $t=5.0$ and $t=10$ are shown.
The first shows the solution obtained by FEM and the second row shows the solution obtained by MMAD.
As expected, due to the dominance of convection, the FEM solution shows an unstable oscillatory behavior which is not observed in the MMAD solution.
At the left boundary $x=0$, the temperature oscillates sinusoidally due to the time-dependent boundary condition. 
This introduces a periodic thermal input into the system.
The flow velocity in the x-direction transports the heat downstream, causing the temperature distribution to shift to the right over time.
Conduction, which is not significant, spreads the heat away from the source, both along the x-direction and across the y-direction. 
Downstream, the temperature oscillations become less pronounced due to the dissipative effect of conduction. 
The solutions highlight the roles played by convection (which carries oscillations) and conduction (which smooths them) over time.
As time elapses, the temperature gradients reduce, leading to a more uniform distribution in the domain.


\section{Summary and outlook}\label{sec:conc}
\noindent
We proposed a novel micromorphic-based artificial diffusion methodology to circumvent the instabilities associated with standard Galerkin finite element approximations of convection-diffusion problems.
Via introducing a new micromorphic-type variable, we add to the formulation an enhancement in the form of terms involving this new variable and its gradient.
The resulting formulation is shown to be stable and convergent also in the pure advection limit.
Comparison with well-established methodologies such as SUPG demonstrates consistent stability and high accuracy of our approach not only for convection-dominated problems but also for transient advection.
Our proposed scheme turns out to be the most accurate methodology when it comes to maintaining the solution characteristics in both convection-dominated and transient advection problems.
Further work includes extension to coupled field problems such as thermofluid flows.


\section*{Acknowledgement}
\noindent
Soheil Firooz and Paul Steinmann gratefully acknowledge the support by the Deutsche Forschungsgemeinschaft (DFG, German Research Foundation) project number 460333672 – CRC 1540 Exploring Brain Mechanics (subproject C01).


\bibliography{library}

\begin{thebibliography}{10}
\expandafter\ifx\csname url\endcsname\relax
  \def\url#1{\texttt{#1}}\fi
\expandafter\ifx\csname urlprefix\endcsname\relax\def\urlprefix{URL }\fi
\expandafter\ifx\csname href\endcsname\relax
  \def\href#1#2{#2} \def\path#1{#1}\fi

\bibitem{Christie1976}
I.~Christie, D.~F. Griffiths, A.~R. Mitchell, O.~C. Zienkiewicz, {Finite
  element methods for second order differential equations with significant
  first derivatives}, International Journal for Numerical Methods in
  Engineering 10 (1976) 1389--1396.

\bibitem{Heinrich1977}
J.~C. Heinrich, P.~S. Huyakorn, O.~C. Zienkiewicz, A.~R. Mitchell, {An
  ‘upwind' finite element scheme for two‐dimensional convective transport
  equation}, International Journal for Numerical Methods in Engineering 11
  (1977) 131--143.

\bibitem{Johnson1984}
C.~Johnson, U.~N{\"{a}}vert, J.~Pitk{\"{a}}ranta, {Finite element methods for
  linear hyperbolic problems}, Computer Methods in Applied Mechanics and
  Engineering 45 (1984) 285--312.

\bibitem{Gresho1981}
P.~M. Gresho, R.~L. Lee, {Don't suppress the wiggles-They're telling you
  something!}, Computers and Fluids 9 (1981) 223--253.

\bibitem{Brooks1982}
A.~N. Brooks, T.~J.~R. Hughes, {Streamline upwind/Petrov-Galerkin formulations
  for convection dominated flows with particular emphasis on the incompressible
  Navier-Stokes equations.}, Computer Methods in Applied Mechanics and
  Engineering 32 (1982) 199--259.

\bibitem{Hughes1978}
T.~J.~R. Hughes, {A simple scheme for developing ‘upwind' finite elements},
  International Journal for Numerical Methods in Engineering 12 (1978)
  1359--1365.

\bibitem{Hughes1982b}
T.~J.~R. Hughes, {A theoretical framework for Petrov-Galerkin methods with
  discontinuous weighting functions: Application to the streamline-upwind
  procedure}, in: R.~H. Gallagher, D.~H. Norrie, J.~T. Oden, O.~C. Zienkiewicz
  (Eds.), Finite element in fluids, John Wiley, 1982, Ch.~3.

\bibitem{Burman2010}
E.~Burman, {Consistent SUPG-method for transient transport problems: Stability
  and convergence}, Computer Methods in Applied Mechanics and Engineering 199
  (2010) 1114--1123.

\bibitem{Hughes1986}
T.~J.~R. Hughes, M.~Mallet, {A new finite element formulation for computational
  fluid dynamics: IV. A discontinuity-capturing operator for multidimensional
  advective-diffusive systems}, Computer Methods in Applied Mechanics and
  Engineering 58 (1986) 329--336.

\bibitem{Hughes1986d}
T.~J.~R. Hughes, M.~Mallet, A.~Mizukami, {A new finite element formulation for
  computational fluid dynamics: II. Beyond SUPG}, Computer Methods in Applied
  Mechanics and Engineering 54 (1986) 341--355.

\bibitem{Knopp2002}
T.~Knopp, G.~Lube, G.~Rapin, {Stabilized finite element methods with shock
  capturing for advection-diffusion problems}, Computer Methods in Applied
  Mechanics and Engineering 191 (2002) 2997--3013.

\bibitem{DeSampaio2001}
P.~A.~B. {De Sampaio}, A.~L. G.~A. Coutinho, {A natural derivation of
  discontinuity capturing operator for convection-diffusion problems}, Computer
  Methods in Applied Mechanics and Engineering 190 (2001) 6291--6308.

\bibitem{Codina1993a}
R.~Codina, {A discontinuity-capturing crosswind-dissipation for the finite
  element solution of the convection-diffusion equation}, Computer Methods in
  Applied Mechanics and Engineering 110 (1993) 325--342.

\bibitem{Hughes1995}
T.~J.~R. Hughes, {Multiscale phenomena: Green's functions, the
  Dirichlet-to-Neumann formulation, subgrid scale models, bubbles and the
  origins of stabilized methods}, Computer Methods in Applied Mechanics and
  Engineering 127 (1995) 387--401.

\bibitem{Hughes1998}
T.~J.~R. Hughes, G.~R. Feij{\'{o}}o, L.~Mazzei, J.~B. Quincy, {The variational
  multiscale method - A paradigm for computational mechanics}, Computer Methods
  in Applied Mechanics and Engineering 166 (1998) 3--24.

\bibitem{Buffa2006}
A.~Buffa, T.~J.~R. Hughes, G.~Sangalli, {Analysis of a multiscale discontinuous
  Galerkin method for convection-diffusion problems}, SIAM Journal on Numerical
  Analysis 44 (2006) 1420--1440.

\bibitem{John2008a}
V.~John, S.~Kaya, {Finite element error analysis of a variational multiscale
  method for the Navier-Stokes equations}, Advances in Computational
  Mathematics 28 (2008) 43--61.

\bibitem{Layton2002}
W.~Layton, {A connection between subgrid scale eddy viscosity and mixed
  methods}, Applied Mathematics and Computation 133 (2002) 147--157.

\bibitem{John2006a}
V.~John, S.~Kaya, W.~Layton, {A two-level variational multiscale method for
  convection-dominated convection-diffusion equations}, Computer Methods in
  Applied Mechanics and Engineering 195 (2006) 4594--4603.

\bibitem{Chen2018a}
Z.~J. Chen, Z.~Y. Li, W.~L. Xie, X.~H. Wu, {A two-level variational multiscale
  meshless local Petrov–Galerkin (VMS-MLPG) method for convection-diffusion
  problems with large Peclet number}, Computers and Fluids 164 (2018) 73--82.

\bibitem{Du2015}
B.~Du, H.~Su, X.~Feng, {Two-level variational multiscale method based on the
  decoupling approach for the natural convection problem}, International
  Communications in Heat and Mass Transfer 61 (2015) 128--139.

\bibitem{Knobloch2009}
P.~Knobloch, G.~Lube, {Local projection stabilization for
  advection-diffusion-reaction problems: One-level vs. two-level approach},
  Applied Numerical Mathematics 59 (2009) 2891--2907.

\bibitem{Cibik2011}
A.~{\c{C}}ibik, S.~Kaya, {A projection-based stabilized finite element method
  for steady-state natural convection problem}, Journal of Mathematical
  Analysis and Applications 381 (2011) 469--484.

\bibitem{Matthies2009}
G.~Matthies, {Local projection stabilisation for higher order discretisations
  of convection-diffusion problems on Shishkin meshes}, Advances in
  Computational Mathematics 30 (2009) 315--337.

\bibitem{Johnson1986}
C.~Johnson, J.~Pitkaranta, {An analysis of the discontinuous Galerkin method
  for a scalar hyperbolic equation}, Mathematics of Computation 46 (1986)
  1--26.

\bibitem{Cockburn1999c}
B.~Cockburn, {Discontinuous Galerkin methods for convection-dominated
  problems}, in: High-order methods for computational physics, Springer, 1999,
  pp. 69----224.

\bibitem{Hughes1989a}
T.~J.~R. Hughes, L.~P. Franca, G.~M. Hulbert, {A new finite element formulation
  for computational fluid dynamics: VIII. The galerkin/least-squares method for
  advective-diffusive equations}, Computer Methods in Applied Mechanics and
  Engineering 73 (1989) 173--189.

\bibitem{Pironneau1992}
O.~Pironneau, J.~Liou, T.~Tezduyar, {Characteristic-Galerkin and Galerkin /
  least-squares space-time formulations for the advection- diffusion equation
  with time-dependent domains}, Computer Methods in Applied Mechanics and
  Engineering 100 (1992) 117--141.

\bibitem{Burman2009}
E.~Burman, M.~A. Fern{\'{a}}ndez, {Finite element methods with symmetric
  stabilization for the transient convection-diffusion-reaction equation},
  Computer Methods in Applied Mechanics and Engineering 198 (2009) 2508--2519.

\bibitem{Guermond1999}
J.~L. Guermond, {Stabilization of Galerkin approximations of transport
  equations by subgrid modeling}, Mathematical Modelling and Numerical Analysis
  33 (1999) 1293--1316.

\bibitem{Codina2000a}
R.~Codina, {Stabilization of incompressibility and convection through
  orthogonal sub-scales in finite element methods}, Computer Methods in Applied
  Mechanics and Engineering 190 (2000) 1579--1599.

\bibitem{Franca1998b}
L.~P. Franca, A.~Nesliturk, M.~Stynes, {On the stability of residual-free
  bubbles for convection-diffusion problems and their approximation by a
  two-level finite element method}, Computer Methods in Applied Mechanics and
  Engineering 166 (1998) 35--49.

\bibitem{Brezzi1999b}
F.~Brezzi, T.~J.~R. Hughes, L.~D. Marini, A.~Russo, E.~Suli, {A priori error
  analysis of residual-free bubbles for advection-diffusion problems}, SIAM
  Journal on Numerical Analysis 36 (1999) 1933--1948.

\bibitem{Donea1984}
J.~Donea, {A Taylor–Galerkin method for convective transport problems},
  International Journal for Numerical Methods in Engineering 20 (1984)
  101--119.

\bibitem{Taylor1973}
C.~Taylor, P.~Hood, {A numerical solution of the Navier-Stokes equations using
  the finite element technique}, Computers and Fluids 1 (1973) 73--100.

\bibitem{Case2011}
M.~A. Case, V.~J. Ervin, A.~Linke, G.~R. Rebholz, {A connection between
  Scott–Vogelius and grad-div stabilized Taylor–Hood FE approximations of
  the Navier–Stokes equations}, SIAM Journal on Mathematical Analysis 49
  (2011) 1461--1481.

\bibitem{Firooz2024}
S.~Firooz, B.~D. Reddy, V.~Zaburdaev, P.~Steinmann, {Mean zero artificial
  diffusion for stable finite element approximation of convection in cellular
  aggregate formation}, Computer Methods in Applied Mechanics and Engineering
  419 (2024) 116649.

\bibitem{Forest2009}
S.~Forest, {Micromorphic approach for gradient elasticity, viscoplasticity, and
  damage}, Journal of Engineering Mechanics 135 (2009) 117--131.

\bibitem{Forest2016}
S.~Forest, {Nonlinear regularization operators as derived from the micromorphic
  approach to gradient elasticity, viscoplasticity and damage}, Proceedings of
  the Royal Society A: Mathematical, Physical and Engineering Sciences 472
  (2016) 20150755.

\bibitem{Forest2020}
S.~Forest, K.~Sab, {Finite-deformation second-order micromorphic theory and its
  relations to strain and stress gradient models}, Mathematics and Mechanics of
  Solids 25 (2020) 1429--1449.

\bibitem{Brooks1981}
A.~N. Brooks, {A Petrov-Galerkin Finite Element Formulation for Convection
  Dominated Flows}, Ph.D. thesis, California Institute of Technology (1981).

\bibitem{Reddy1998}
B.~D. Reddy, {Introductory Functional Analysis: with Applications to Boundary
  Value Problems and Finite Elements}, Springer Science and Business Media,
  1998.

\end{thebibliography}


\end{document}